\documentclass[preprint,3p,sort&compress,final,times]{elsarticle}
\usepackage{amsmath}
\usepackage{amssymb}
\usepackage{mathtools}
\usepackage{placeins}
\usepackage[usenames,dvipsnames]{color}

\usepackage{cases}

\usepackage{lineno}

\newcommand{\ip}[2]{\langle #2, #1\rangle_{\Omega}}
\renewcommand{\d}{\mathrm{d}}

\newcommand{\figref}[1]{Figure~\ref{#1}} 

\newcommand{\ederiv}{\mathrm{d}}
\newcommand{\vfield}[1]{\vec{#1}}             
\newcommand{\innerproduct}[2]{\langle #1,#2\rangle_{\Omega}}
\newcommand{\gradientperp}{\nabla^{\perp}}      
\newcommand{\ddt}[2][]{\frac{\ederiv^{#1} #2}{\ederiv t^{#1}}}     
\newcommand{\pddt}[2][]{\frac{\partial^{#1} #2}{\partial t^{#1}}}   
\newcommand{\transpose}{\intercal}                      
\newcommand{\valgebra}[1]{\boldsymbol{#1}}                             
\newcommand{\valgebracomponents}[1]{#1}                             
\newcommand{\moperator}[1]{\boldsymbol{\mathsf{#1}}}            
\newcommand{\moperatorcomponents}[1]{\mathsf{#1}}                                  

\newcommand{\firstRev}[1]{#1}
\newcommand{\secondRev}[1]{#1}

\newdefinition{remark}{Remark}

\journal{Journal of Computational Physics}
\begin{document}
\begin{frontmatter}

\title{Discrete conservation properties for shallow water flows using mixed mimetic spectral elements}
\author[CCS]{D.~Lee\corref{cor}}
\ead{drlee@lanl.gov}
\author[EUT]{A. Palha}
\author[TUD]{M. Gerritsma}

\address[CCS]{Computer, Computational and Statistical Sciences, Los Alamos National Laboratory, Los Alamos, NM 87545, USA}
\address[EUT]{Eindhoven University of Technology, Department of Mechanical Engineering, P.O. Box 513, 5600 MB Eindhoven, The Netherlands}
\address[TUD]{Delft University of Technology, Faculty of Aerospace Engineering, P.O. Box 5058, 2600 GB Delft, The Netherlands}
\cortext[cor]{Corresponding author. Tel. +1 505 665 7286.}

\begin{abstract}
A mixed mimetic spectral element method is applied to solve the rotating shallow water equations. The mixed method uses
the recently developed spectral element histopolation functions, which exactly satisfy the fundamental theorem of calculus with
respect to the standard Lagrange basis functions in one dimension. These are used to construct tensor product solution spaces
which satisfy the generalized Stokes theorem, as well as the annihilation of the gradient operator by the curl and the curl by the divergence. 
This allows for the exact conservation of first order moments (mass, vorticity), as well as \secondRev{higher } moments (energy, potential enstrophy),
subject to the truncation error of the time stepping scheme. The continuity equation is solved in the strong form, such that
mass conservation holds point wise, while the momentum equation is solved in the weak form such that vorticity is globally
conserved. While mass, vorticity and energy conservation hold for any quadrature rule, potential
enstrophy conservation is dependent on exact spatial integration.
The method possesses a weak form statement of geostrophic balance due to the compatible nature of the solution spaces
and arbitrarily high order spatial error convergence.

\end{abstract}

\begin{keyword}
Mimetic\sep
Spectral elements\sep
High order\sep
Shallow water\sep
Energy and potential enstrophy conservation
\end{keyword}

\end{frontmatter}

\section{Introduction}

In recent years there has been much interest in the use of finite element methods for the
development of geophysical fluid solvers. This is in large part due to the recognition of the
importance of conservation over long time integrations as a means of mitigating against
both numerical instabilities and biases in the solution \cite{Thuburn08}, and the capacity
of finite elements to satisfy the conservation of various moments via the use of compatible
or mimetic finite element spaces \cite{PG17}. Various finite element spaces have been explored 
for their suitability for modelling geophysical flows including Raviart-Thomas, Brezzi-Douglas-Marini
and Brezzi-Douglas-Fortin-Marini elements \cite{CS12,MC14,NSC16}. When used in a
sequence of element types such as the standard continuous and discontinuous Galerkin elements,
these can be shown to exactly satisfy the Kelvin-Stokes and Gauss-divergence theorems
when applying the curl and divergence operators respectively in the weak form, as well as the
annihilation of the gradient operator by the curl. Satisfying these properties exactly in the
discrete form is necessary in order to conserve both \secondRev{first order and higher order } 
moments for the shallow water equations in rotational form, as first presented for a C-grid finite 
difference scheme \cite{AL81}, and later formalized via the derivation of the finite difference 
operators from Hamiltonian methods \cite{Salmon04,Salmon07}.

The application of mimetic finite elements for shallow water flows has also been formalized in
the language of exterior calculus in order to generalize the expression of their conservation
properties \cite{CT14}. Mimetic properties have also been demonstrated for the standard A-grid
spectral element method on cubed sphere geometries via careful use of covariant and contra-variant
transformations in order to evaluate the curl and divergence operators respectively so as to
exactly satisfy the Kelvin-Stokes and Gauss-divergence theorems \cite{TF10}.

The present study explores the use of spectral elements using mixed basis functions in order
to preserve mimetic properties for geophysical flows. These
recently developed methods invoke the use of \emph{histopolation} functions \cite{Gerritsma11}, which are
defined such that they exactly satisfy the fundamental theorem of calculus with respect to the
nodal Lagrange basis functions from which they are derived. In the language of exterior
calculus these edge functions may then be regarded as defining the space of 1-forms (with the
Lagrange polynomials defining the space of 0-forms). By taking tensor product combinations
of these basis functions and their associated edge functions higher dimensional differential $k$-forms may
also be constructed for which the Kelvin-Stokes and Gauss-divergence theorems are satisfied \cite{KG13}.

As for other compatible tensor product finite element methods, the mixed mimetic spectral elements
provide an effective means of preserving many of the properties of a geophysical fluid over long time
integrations. These include:

\begin{itemize}
\item \emph{Conservation of mass, vorticity, total energy and potential enstrophy}\\
Conservation of mass ensures that the solution does not drift over time and develop
large biases. Vorticity is an important moment to conserve since the large scale circulation
of the atmosphere and oceans at mid latitudes is dominated by a quasi-balance between rotation and
pressure gradients. Conservation of energy ensures that solutions remain bounded
and numerical instabilities do not grow exponentially. The importance of potential enstrophy
conservation is less immediately obvious meanwhile since for two dimensional turbulent fluids this
cascades to small scales where some dissipation mechanism is necessary \cite{Vallis06}.
\item \emph{Stationary geostrophic modes}\\
The large scale circulations of the atmosphere and ocean are dominated by the slow evolution
of modes balanced between Coriolis forces and pressure gradients. If this balance cannot be
exactly replicated in a numerical model then the process of adjustment will result in the
radiation of fast gravity waves that can contaminate the solution \cite{TRSK09}.
\item \emph{High order spatial error convergence}\\
The spectral element method defines the nodes of the basis functions to cluster towards element boundaries
such that spurious oscillations due to spectral ringing are avoided, allowing convergence of errors at
arbitrarily high order.
\end{itemize}

The rest of this paper proceeds as follows: In the following section the shallow water equations will be briefly
introduced in the continuous form. Section 3 discusses the mixed mimetic spectral element method,
as introduced by previous authors \cite{Gerritsma11,KG13,KPG11,Hiemstra14}, and its application to
the rotating shallow water equations. Section 4 explores the conservation properties of the discrete system.
In Section 5 some results are presented, demonstrating
the error convergence, conservation and balance properties of the method. Section 6 details some conclusions
regarding the suitability of the method for geophysical flows, its advantages and limitations, as well as some
future work we intend to pursue on the topic.

\section{The 2D shallow water equations}
	The two dimensional shallow water equations present an excellent testing ground for primitive equation atmospheric and
oceanic models since they exhibit many of the same features, including slowly evolving Rossby waves
and fast gravity waves, nonlinear cascades of \secondRev{higher } moments (kinetic energy and potential enstrophy),
and the conservation of various moments (mass, total energy, vorticity, and an infinite
number of rotational moments, of which potential enstrophy is the first).

	Let $\Omega = [x_{\mathrm{min}},x_{\mathrm{max}}]\times[y_{\mathrm{min}},y_{\mathrm{max}}]\subset \mathbb{R}^{2}$ such that
$x_{\mathrm{min}} < x_{\mathrm{max}}$ and $y_{\mathrm{min}} < y_{\mathrm{max}}$ be a doubly periodic domain, and let $t_{F}>0$.
The rotational form of the shallow water problem, \cite{AL81}, consists in finding the prognostic variables velocity,
$\vec{u}:\Omega\times(0,t_{F}]\mapsto\mathbb{R}^{2}$, and fluid depth, $h:\Omega\times(0,t_{F}]\mapsto\mathbb{R}$, such that
\footnote{The shallow water equations are formulated in 2D, such that $q \times \vec{F}:= q\vec{e}_z \times \vec{F} = q \vec{F}^{\perp}$.}
\begin{subnumcases}{\label{eq:shallow_water_equations_all}}
	\frac{\partial\vec u}{\partial t} + q\times\vec F + \nabla(K + gh) = 0\,, & in $\Omega\times (0,t_{F}]\,,$  \label{mom_cont}\\
	\frac{\partial h}{\partial t} + \nabla\cdot\vec F = 0\,, & in $\Omega\times (0,t_{F}]\,,$ \label{mas_cont}
\end{subnumcases}
The diagnostic variables are potential vorticity $q:\Omega\times(0,t_{F}]\mapsto\mathbb{R}$, mass flux $\vec{F}:\Omega\times(0,t_{F}]\mapsto\mathbb{R}^{2}$,
and kinetic energy per unit mass $K:\Omega\times(0,t_{F}]\mapsto\mathbb{R}$, defined as
\begin{subnumcases}{\label{eq:prognostic_variables}}
	q := \frac{\nabla\times\vec u + f}{h}\,, &  \label{eq:definition_potential_vorticity} \\
	\vec{F} := h\vec{u}\,, &  \label{eq:definition_h_flux} \\
	K := \frac{1}{2}\vec u\cdot\vec u\,, &  \label{eq:definition_kinetic_energy}
\end{subnumcases}
where $f:\Omega\mapsto\mathbb{R}$ is the Coriolis term. Note that we assume that this Coriolis term does not explicitly depend on time.

The shallow water equations conserve the following invariants, \cite{AL81}:
\begin{itemize}
\item Volume: integrating (\ref{mas_cont}) over the domain $\Omega$ and assuming periodic
boundary conditions gives the result that
\begin{equation}
	\frac{\mathrm{d}}{\mathrm{d} t}\int_{\Omega}h \,\d\Omega = 0\,.
\end{equation}
If we assume constant density, then this also gives mass conservation.

\item Vorticity: taking the curl of (\ref{mom_cont}) leads to a conservation equation for vorticity
$\omega = \nabla\times\vec u:\Omega\times(0,t_{F}]\mapsto\mathbb{R}$, of the form
\begin{equation}
\frac{\partial\omega}{\partial t} + \nabla\cdot(\vec u\,(\omega+f))=0\;.
\end{equation}
Integrating over the domain (and assuming periodic boundary conditions) gives the conservation
of vorticity as
\begin{equation}
	\frac{\mathrm{d}}{\mathrm{d} t}\int_{\Omega}\omega \,\d\Omega = 0\,.
\end{equation}

\item Energy: computing the inner product of (\ref{mom_cont}) with $\vec F$ and multiplying (\ref{mas_cont}) by $gh$ gives
\begin{align}
\frac{\partial hK}{\partial t} + \nabla\cdot(\vec FK) + \vec F\cdot\nabla(gh) &= 0 \;,\\
\frac{\partial(\frac{1}{2}gh^2)}{\partial t} + \nabla\cdot(gh\vec F) - \vec F\cdot\nabla(gh) &= 0 \;.
\end{align}
Combining these and integrating over $\Omega$ gives energy conservation as
\begin{equation}
	\frac{\mathrm{d}}{\mathrm{d} t}\int_{\Omega}\left \{hK + \frac{1}{2}gh^2 \right \}\,\d\Omega = 0\,.
\end{equation}

\item Potential enstrophy: expressing the vorticity as $\omega = hq-f$
    within the
vorticity advection equation and subtracting $q$ times (\ref{mas_cont}) gives an advection
equation for the potential vorticity as

\begin{equation}
\frac{\partial q}{\partial t} + \vec u\cdot\nabla q = 0.
\end{equation}
Multiplying this by $hq$ and (\ref{mas_cont}) by $\frac{1}{2}q^2$ and adding gives

\begin{equation}
\frac{\partial\frac{1}{2}hq^2}{\partial t} + \nabla\cdot\Big(\frac{1}{2}\vec Fq^2\Big) = 0.
\end{equation}
Integrating over the domain then gives potential enstrophy conservation as
\begin{equation}
	\frac{\mathrm{d}}{\mathrm{d} t}\int_{\Omega}\frac{1}{2}hq^2\,\d\Omega = 0\,.
\end{equation}
\end{itemize}
For a detailed derivation of these conservation laws the reader is referred to previous studies
\cite{MC14,AL81}.

One of the central objectives of this paper is to show that these properties can be satisfied at
the discrete level within the mimetic spectral element discretization framework. Ensuring the
conservation of invariants in discrete form helps to mitigate against biases and instabilities
in the solution of the original system \cite{Thuburn08}. For this reason, it is considered important
to satisfy these conservation properties at the discrete level in geophysical flow solvers,
especially when long time simulation is the goal.

\section{Spatial discretization}
	In this work we set out to construct a mimetic spectral finite element discretization for the shallow water equations as given by
\eqref{eq:shallow_water_equations_all}, together with the diagnostic equations \eqref{eq:prognostic_variables}. In particular we use a
mixed finite element formulation, for more details on mixed finite elements see for example \cite{brezzi1991mixed,BoffiMixedFiniteElements2013}.
	
\firstRev{The mimetic finite element discretization presented here differs from other finite element methods in the degrees 
of freedom. Here the degrees of freedom, the expansion coefficients, represent integral values. For nodal basis functions, the 
degrees of freedom represent the values at points, see \eqref{eq::2d_nodal_polynomials_properties}, the degrees of freedom 
for vectors will be the (two-dimensional) fluxes over line segments, see \eqref{eq::edge_polynomials_properties_xi}, 
\eqref{eq::edge_polynomials_properties_eta} in the section. The degrees of freedom for densities will be their integrated 
values over 2D volumes, see \eqref{eq::volume_polynomials_properties}. The main motivation for the use of these degrees of 
freedom are: 1. that integral values are well defined on non-orthogonal grids and 2. that the optimal order of convergence 
is displayed on both affine and non-affine meshes. Alternative formulations may encounter loss of optimality, \cite{ArnoldBoffiBonizzoni,ArnoldAwanou,ArnoldBoffiFalk}. 
An important feature of this sequence of basis functions is that the derivative can be applied directly to the degrees of 
freedom by multiplying the degrees of freedom by an appropriate incidence matrix
\cite{bossavit_japan_computational_1,bossavit_japan_computational_2,bossavit_japan_computational_3,bossavit_japan_computational_4,bossavit_japan_computational_5},
as shown in \eqref{eq:hilbert_subcomplex_basis_W} and \eqref{eq:hilbert_subcomplex_basis_U} in this section.
The incidence matrices do not depend on the polynomials degree and the shape of the grid; they represent the topological 
derivative. In Section~4 we will prove conservation properties in terms of these topological structures. Since the proofs 
are based on metric-free concepts, this ensures that these properties will also hold on highly deformed grids.}

	\subsection{Weak formulation}
		The first step to construct this discretization is the weak form of \eqref{eq:shallow_water_equations_all} and \eqref{eq:prognostic_variables}.
		In this work, as usual, $\ip{\cdot}{\cdot}$ represents the $L^{2}$ inner product
		\begin{equation}
			\ip{g}{f} := \int_{\Omega} f\cdot g\,\mathrm{d}\Omega\,.
		\end{equation}
		The weak formulation reads: for any time $t\in (0,t_{F}]$ and for a given Coriolis term $f\in L^2(\Omega)$, find
		$\vec{u},\vec{F}\in H(\mathrm{div},\Omega)$, $h,K\in L^{2}(\Omega)$, and $q\in H(\mathrm{rot},\Omega)$
\footnote{For scalar variable $\psi$ the rot operator $\nabla^{\perp} := \vec e_z\times\nabla$ is defined as
\[ \nabla^{\perp}\psi = \left ( \begin{array}{c} -\partial \psi/\partial y \\
\partial \psi/\partial x
\end{array} \right ) \;.\]} such that
	 \begin{subnumcases}{\label{eq:shallow_water_continuous_weak_form}}
		\ip{\frac{\partial\vec{u}}{\partial t}}{\vec{\nu}} + \ip{q\times\vec F}{\vec{\nu}} - \ip{K + gh}{\nabla\cdot\vec{\nu}} = 0\,, &
			$\forall \vec{\nu}\in H(\mathrm{div},\Omega)$  \label{mom_cont_weak} \\
		\ip{\frac{\partial h}{\partial t}}{\sigma} + \ip{\nabla\cdot\vec F}{\sigma} = 0\,, & $\forall \sigma\in L^{2}(\Omega)$, \label{mas_cont_weak} \\
		\ip{hq}{\zeta} = -\ip{\vec u}{\nabla^{\perp}\zeta} + \ip{f}{\zeta}\,, & $\forall \zeta\in H(\mathrm{rot},\Omega)$, \label{eq:definition_potential_vorticity_weak} \\
		\ip{\vec{F}}{\vec{\varphi}} = \ip{h\vec{u}}{\vec{\varphi}}\,, &  $\forall \vec{\varphi}\in H(\mathrm{div},\Omega)$, \label{eq:definition_h_flux_weak} \\
		\ip{K}{\kappa} = \frac{1}{2}\ip{\vec{u}\cdot\vec{u}}{\kappa}\,, & $\forall \kappa\in L^{2}(\Omega)$,  \label{eq:definition_kinetic_energy_weak}
	\end{subnumcases}
	where we have used integration by parts and the periodic boundary conditions to obtain the identities \cite{MC14}
	 \begin{subnumcases}{\label{eq:adjoints}}
	\ip{\nabla(K + gh)}{\vec{\nu}} = -\ip{K + gh}{\nabla\cdot\vec{\nu}},\label{adjoint_divergence}\\
	\ip{\nabla\times\vec{u}}{\varsigma} = -\ip{\vec{u}}{\nabla^{\perp}\varsigma}\label{adjoint_curl}.
	\end{subnumcases}
	These two relations show that minus the gradient is the Hilbert adjoint of the divergence operator and minus the curl is the
	Hilbert adjoint of the rot, provided the boundary integrals vanish. The space $L^{2}(\Omega)$ corresponds to square
	integrable functions and the spaces $H(\mathrm{div},\Omega)$ and $H(\mathrm{rot},\Omega)$ contain square integrable
	functions whose divergence and rot are also square integrable.
	
	\subsection{Finite dimensional mimetic function spaces} \label{sec::the_finite_dimensional_basis_functions}
		The second step to construct this discretization is the definition of the spatial conforming function spaces,
		where we will seek the discrete solutions for velocity $\vec{u}_{h}$, fluid depth $h_{h}$, mass flux $\vec{F}_{h}$,
		kinetic energy $K_{h}$ and potential vorticity $q_{h}$:
		\begin{equation}
			q_{h}\in W_{h}\subset H(\mathrm{rot},\Omega),\quad \vec{u}_{h},\vec{F}_{h}\in U_{h}\subset H(\mathrm{div},\Omega), \quad
			\mathrm{and} \quad h_{h},K_{h}\in Q_{h}\subset L^{2}(\Omega)\,.
		\end{equation}
		
		The choice of finite dimensional function spaces determines the properties of the discretization, see for example \cite{arnold2006finite,arnold2010finite}
		for a more general discussion, and \cite{PG17} for a recent discussion focussed on the 2D Navier-Stokes equations.
		
		Therefore, the finite dimensional function spaces used in this work are such that when combined form a Hilbert subcomplex
		\begin{equation}
			\mathbb{R} \longrightarrow W_{h} \stackrel{\nabla^{\perp}}{\longrightarrow} U_{h} \stackrel{\nabla\cdot}{\longrightarrow} Q_{h}
			\longrightarrow 0\,.\label{eq:hilbert_subcomplex}
		\end{equation}
		The meaning of this Hilbert subcomplex is that
		\begin{equation}
			 \{\nabla^{\perp} \omega_{h}\,|\,\omega_{h}\in W_{h}\} \subset U_{h} \quad \mathrm{and}\quad \{\nabla\cdot\vec{u}_{h}\,|\,\vec{u}_{h}\in U_{h}\} \subseteq Q_{h}\,.
		\end{equation}
		In other words, the rot operator must map $W_{h}$ into $U_{h}$ and the div operator must map $U_{h}$ onto $Q_{h}$.
		
		This discrete subcomplex mimics the 2D Hilbert complex associated to the continuous functional spaces
		\begin{equation}
			\mathbb{R} \longrightarrow H(\mathrm{rot},\Omega) \stackrel{\nabla^{\perp}}{\longrightarrow} H(\mathrm{div},\Omega)
			\stackrel{\nabla\cdot}{\longrightarrow} L^{2}(\Omega) \longrightarrow 0\,.\label{eq:hilbert_complex}
		\end{equation}
		The Hilbert complex is an important structure that is intimately connected to the de Rham complex of differential forms.
		Therefore, the construction of a discrete subcomplex is an important requirement to obtain a stable and accurate finite
		element discretization, see for example
		 \cite{arnold2010finite,Palha2014,{bossavit_japan_computational_1,bossavit_japan_computational_2,bossavit_japan_computational_3,bossavit_japan_computational_4,bossavit_japan_computational_5}}
		for a detailed discussion.
		
		Each of these finite dimensional function spaces $W_{h}$, $U_{h}$, and $Q_{h}$ has an associated finite set of basis functions
		$\epsilon_{i}^{\,W}, \vec{\epsilon}_{i}^{\,U}, \epsilon_{i}^{Q}$, such that
		\begin{equation}
			W_{h} = \mathrm{span}\{\epsilon_{1}^{\,W},\dots,\epsilon_{d_{W}}^{\,W}\}, \quad U_{h} = \mathrm{span}\{\vec{\epsilon}_{1}^{\,U},\dots,\vec{\epsilon}_{d_{U}}^{\,U}\}, \quad\mathrm{and}\quad Q_{h} = \mathrm{span}\{\epsilon_{1}^{Q},\dots,\epsilon_{d_{Q}}^{Q}\}\,,
		\end{equation}
		where $d_{W}$, $d_{U}$, and $d_{Q}$ denote the dimension of the discrete function spaces and therefore correspond to the number of
		degrees of freedom associated to each of the unknowns. As a consequence, the approximate solutions for vorticity, potential vorticity,
		Coriolis, velocity, mass flux, fluid depth, and kinetic energy can be expressed as a linear combination of these basis functions
			 \begin{subnumcases}{\label{eq:physical_quantities_expansion}}
				q_{h} := \sum_{i=1}^{d_{W}} q_{i}\,\epsilon_{i}^{\,W},\quad f_h := \sum_{i=1}^{d_W}f_i\epsilon_i^{\,W}, & \label{eq:omega_q_expansion} \\
				\vec{u}_{h} := \sum_{i=1}^{d_{U}} u_{i}\,\vec{\epsilon}_{i}^{\,U},\quad \vec{F}_{h} := \sum_{i=1}^{d_{U}} F_{i}\,\vec{\epsilon}_{i}^{\,U}, & \label{eq:u_F_expansion} \\
				h_{h} := \sum_{i=1}^{d_{Q}} h_{i}\, \epsilon_{i}^{Q},\quad K_{h} := \sum_{i=1}^{d_{Q}} K_{i}\,\epsilon_{i}^{Q}, & \label{eq:h_K_expansion}
			\end{subnumcases}
		with $q_{i}$, $f_i$, $u_{i}$, $F_{i}$, $h_{i}$, $K_{i}$, the degrees of freedom for vorticity, potential vorticity,
		Coriolis, velocity, mass flux, fluid depth, and kinetic energy, respectively. Since the shallow water equations form a time dependent
		set of equations, these coefficients will be time dependent.

\subsubsection{One dimensional nodal and histopolant polynomials}
To define the two-dimensional basis functions $\epsilon_{i}^{\,W}$, $\vec{\epsilon}_{i}^{\,U}$, and $\epsilon_{i}^{Q}$, we first introduce two types
of one-dimensional polynomials: one associated with nodal interpolation, and the other with integral interpolation (histopolation)
\footnote{For an extensive discussion of integral interpolation (histopolation) see \cite{robidoux-polynomial,Gerritsma11}.}. Subsequently, these two
types of polynomials will be combined to generate the two-dimensional polynomial basis functions used to discretize the physical quantities that appear 
in this problem.
			
Consider the canonical interval $I=[-1,1]\subset\mathbb{R}$ and the Legendre polynomials, $L_{p}(\xi)$ of degree $p$ with $\xi\in I$. The $p+1$ roots,
$\xi_{i}$, of the polynomial $(1-\xi^{2})\frac{\mathrm{d}L_{p}}{\mathrm{d}\xi}$ are called Gauss-Lobatto-Legendre (GLL) nodes and satisfy
$-1 = \xi_{0} < \xi_{1} < \dots < \xi_{p-1} < \xi_{p} = 1$. Let $l^{p}_{i}(\xi)$ be the Lagrange polynomial of degree $p$ through the GLL nodes, such that
			\begin{equation}
				l^{p}_{i}(\xi_{j}) :=
				\begin{cases}
					1 & \mbox{if } i = j \\
					& \\
					0 & \mbox{if } i \neq j
				\end{cases}\,, \quad i,j = 0,\dots,p\,. \label{eq::nodal_basis_polynomials}
			\end{equation}
			The explicit form of these Lagrange polynomials is given by
			\begin{equation}
				l^{p}_{i}(\xi) = \prod_{\substack{k=0\\k\neq i}}^{p}\frac{\xi-\xi_{k}}{\xi_{i}-\xi_{k}}\,. \label{eq::lagrange_interpolants}
			\end{equation}

Let $q_h(\xi)$ be a polynomial of degree $p$ defined on $I=[-1,1]$ and $q_{i} = q_h(\xi_{i})$, then the expansion of $q_h(\xi)$ in terms of Lagrange 
polynomials is given by
			\begin{equation}
				q_h(\xi) := \sum_{i=0}^{p}q_{i}l^{p}_{i}(\xi)\,. \label{eq::nodal_polynomial_expansion}
			\end{equation}

Because the expansion coefficients in \eqref{eq::nodal_polynomial_expansion} are given by the value of $q_h$ in the nodes $\xi_i$,
we refer to this interpolation as a {\em nodal interpolation} and we will denote the Lagrange polynomials in \eqref{eq::lagrange_interpolants} by \emph{nodal polynomials}.
			
\begin{figure}[!ht]
				\begin{center}
					 \includegraphics[width=0.5\textwidth]{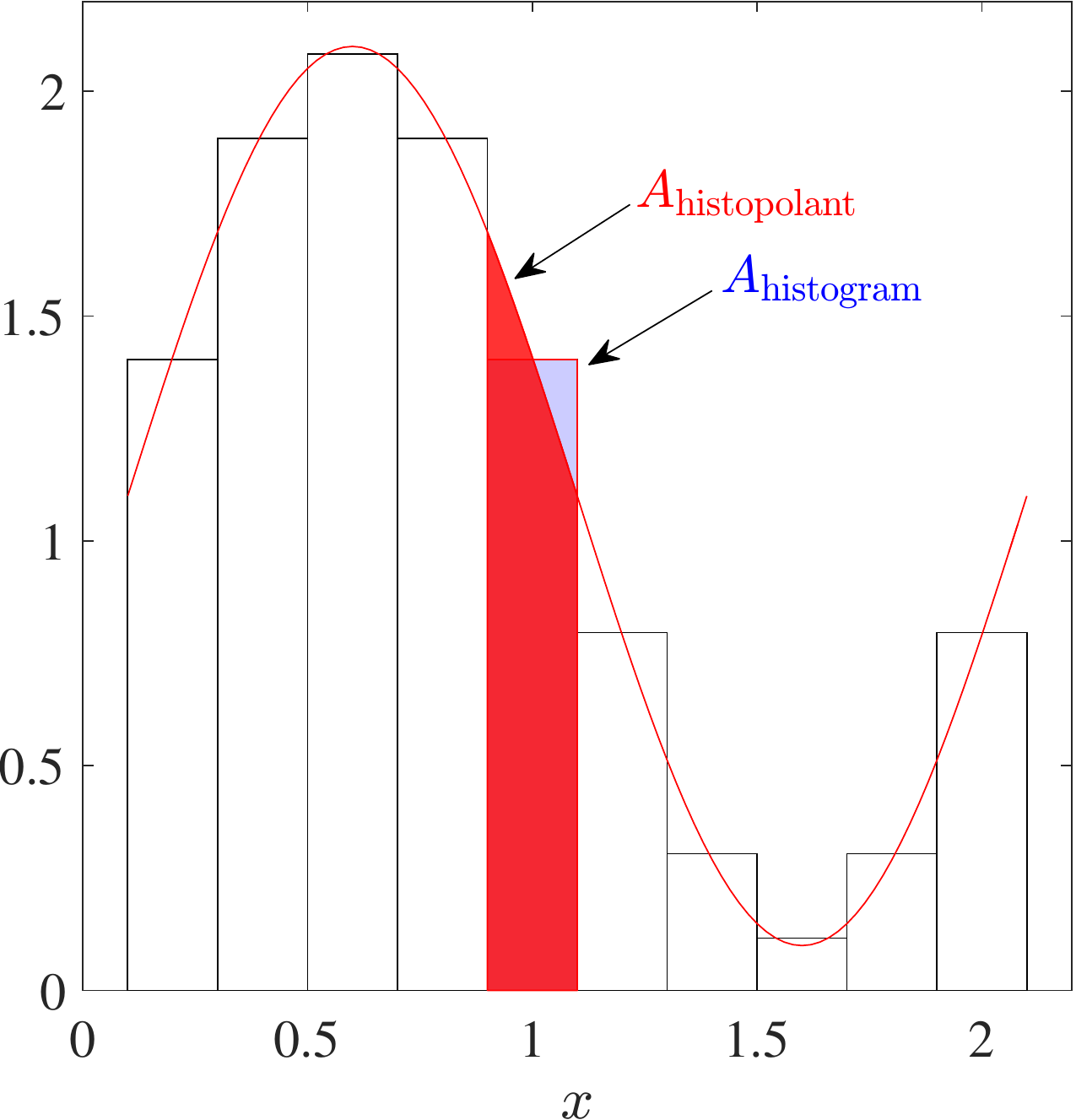}
					\caption{Histogram and an example of a histopolant (red curve). By definition, the integral of the histopolant over each cell (or bin) $A_{\mathrm{histopolant}}$ is equal to the area of the corresponding bar of the histogram $A_{\mathrm{histogram}}$.}
					\label{fig::example_histopolation}
				\end{center}
			\end{figure}

			Before introducing the second set of basis polynomials that will be used in this work it is important to introduce the reader to the concept of \emph{histopolant}. Given a histogram, i.e. a piece-wise constant function, a histopolant is a smooth function whose integrals over the cells (or bins) of the histogram are equal to the area of the corresponding bars of the histogram, see \figref{fig::example_histopolation}. If the histopolant is a polynomial we say that it is a \emph{polynomial histopolant}. In the same way as a polynomial interpolant that passes exactly through $p+1$ points has degree $p$, a polynomial that exactly histopolates a histogram with $p+1$ bins has polynomial degree $p$. Consider now a function $g(x)$ and its associated integrals over a set of cells $[a_{j-1},a_{j}]$, $g_{j} = \int_{a_{j-1}}^{a_{j}} g(x)\mathrm{d}x$, with $ a_{0}, < \dots < a_{j} < \dots < a_{p}$. The set of integral values $g_{j}$ and cells $[a_{j-1},a_{j}]$ can be seen as a histogram. As mentioned before, it is possible to construct a histopolant of this histogram. This histopolant will be an approximating function of $g$ that has the particular property of having the same integral over the cells $[a_{j-1},a_{j}]$ as the original $g$. Just like a nodal interpolation exactly reconstructs the original function at the interpolating points, a histopolant exactly reconstructs the integral of the original function over the cells.
			
			Using the nodal polynomials we can define another set of basis polynomials, $e^{p}_{i}(\xi)$, as
			\begin{equation}
				e^{p}_{i}(\xi) := - \sum_{k=0}^{i-1}\frac{\mathrm{d}l^{p}_{k}(\xi)}{\mathrm{d}\xi}\,, \qquad i=1,\dots,p\,. \label{eq::histopolant_polynomials_definition}
			\end{equation}
			These polynomials $e^{p}_{i}(\xi)$ have polynomial degree $p-1$ and satisfy,
			\begin{equation}
				 \int_{\xi_{j-1}}^{\xi_{j}}e^{p}_{i}(\xi)\,\mathrm{d}\xi =
				\begin{cases}
					1 & \mbox{if } i=j \\
					& \\
					0 & \mbox{if } i\neq j
				\end{cases}\,, \quad i,j = 1,\dots,p\,. \label{eq::histopolant_polynomials_properties}
			\end{equation}
			The proof that the polynomials $e^{p}_{i}(\xi)$ have degree $p-1$ follows directly from the fact that their definition \eqref{eq::histopolant_polynomials_definition} involves a linear combination of the derivative of polynomials of degree $p$. The proof of \eqref{eq::histopolant_polynomials_properties} results from the properties of $l_{k}^{p}(\xi)$. Using  \eqref{eq::histopolant_polynomials_definition} the integral of $e^{p}_{i}(\xi)$ becomes
			\[
				 \int_{\xi_{j-1}}^{\xi_{j}}e^{p}_{i}(\xi)\,\mathrm{d}\xi =  - \int_{\xi_{j-1}}^{\xi_{j}}\sum_{k=0}^{i-1}\frac{\mathrm{d}l^{p}_{k}(\xi)}{\mathrm{d}\xi} = - \sum_{k=0}^{i-1}\int_{\xi_{j-1}}^{\xi_{j}}\frac{\mathrm{d}l^{p}_{k}(\xi)}{\mathrm{d}\xi} = - \sum_{k=0}^{i-1} \left(l^{p}_{k}(\xi_{j}) - l^{p}_{k}(\xi_{j-1})\right) = - \sum_{k=0}^{i-1} \left(\delta_{k,j} -\delta_{k,j-1}\right)\,,
			\]
			where $\delta_{i,j}$ is the Kronecker delta. It is straightforward to see that
			\[
				- \sum_{k=0}^{i-1} \left(\delta_{k,j} -\delta_{k,j-1}\right) = \begin{cases}
					1 & \mbox{if } i=j \\
					& \\
					0 & \mbox{if } i\neq j
				\end{cases}\,, \quad i,j = 1,\dots,p\,.
			\]
			For more details see \cite{robidoux-polynomial,Gerritsma11}.
			
			Let $g_h(\xi)$ be a polynomial of degree $(p-1)$ defined on $I=[-1,1]$ and $g_{i} = \int_{\xi_{i-1}}^{\xi_{i}}g_h(\xi)\,\mathrm{d}\xi$, then its expansion in terms of the polynomials $e_i^p(\xi)$ is given by
			\begin{equation}
				g_{h}(\xi) = \sum_{i=1}^{p}g_{i}e^{p}_{i}(\xi)\;.
\label{eq:1D_expansion_in_edge_polynomials}
			\end{equation}
Because the expansion coefficients in \eqref{eq:1D_expansion_in_edge_polynomials} are the integral values of $g_h(\xi)$, we denote the polynomials in
\eqref{eq::histopolant_polynomials_definition} by \emph{histopolant polynomials}\footnote{In earlier work we referred to these functions as
{\em edge functions}.} and refer to \eqref{eq:1D_expansion_in_edge_polynomials} as {\em histopolation}.
			
\begin{figure}
				\begin{center}
					 \includegraphics[width=0.4\textwidth]{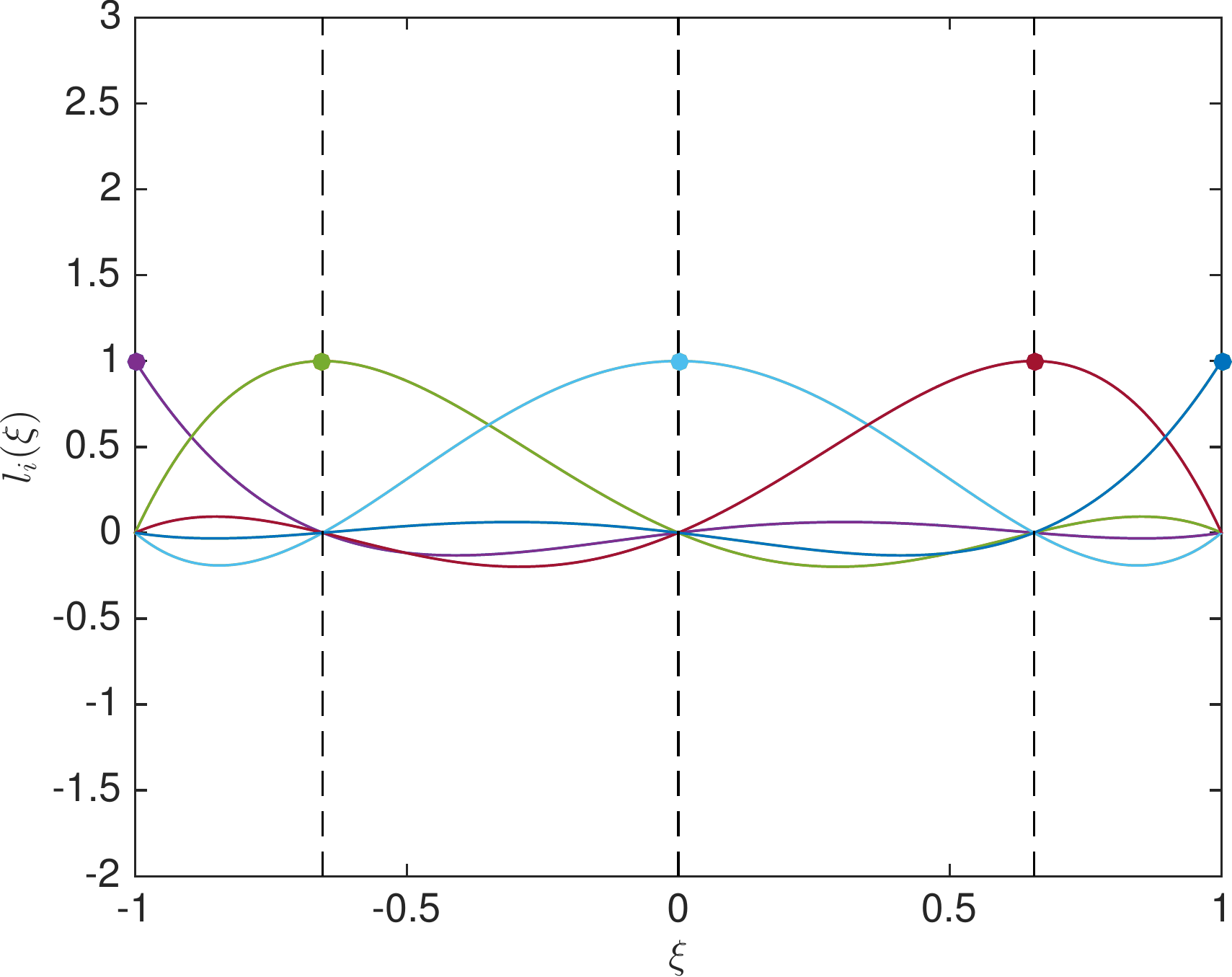} \hspace{1cm}
					 \includegraphics[width=0.4\textwidth]{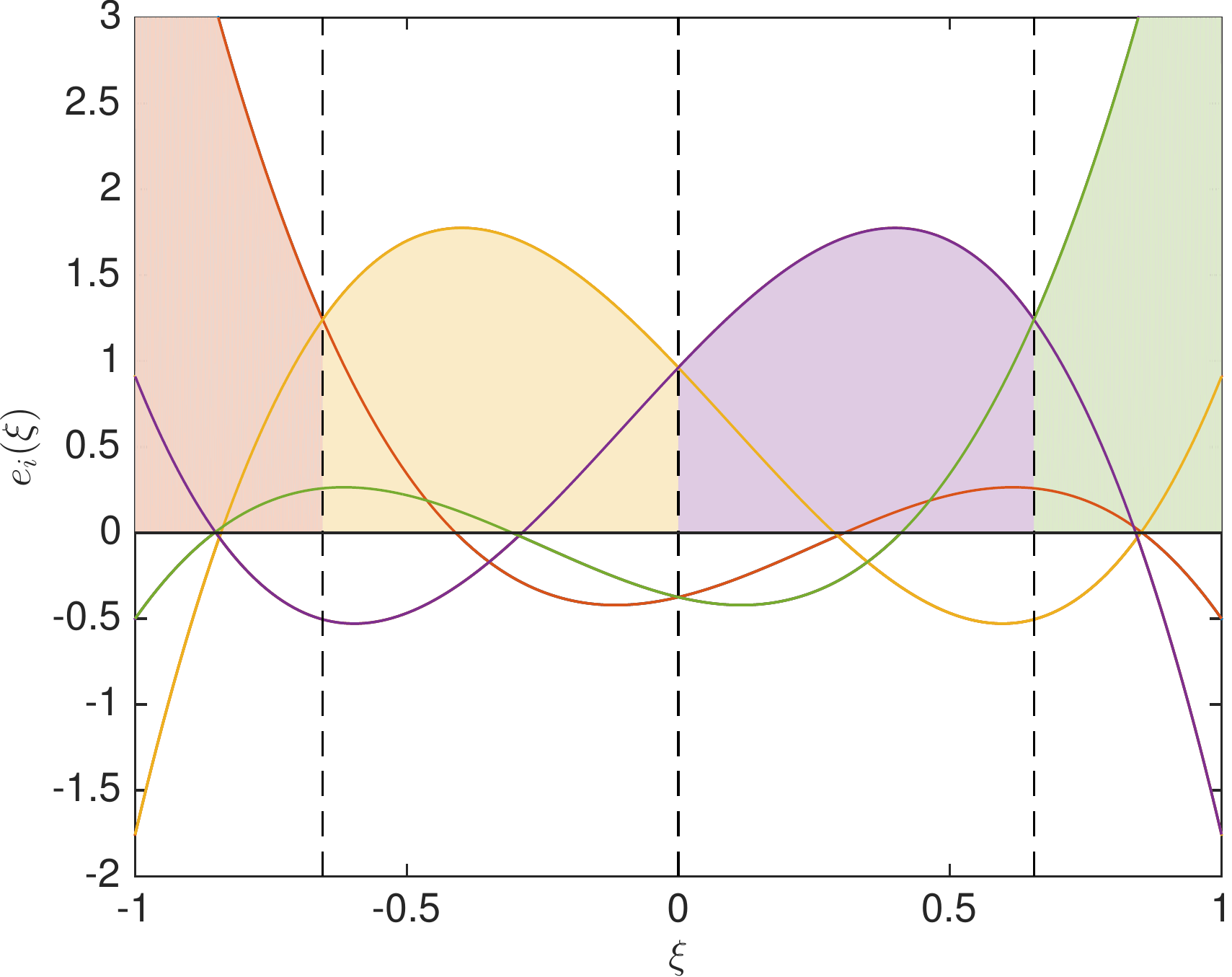}
					\caption{Basis polynomials associated to $p=4$. Left: nodal polynomials, the value of the basis polynomial at the corresponding node is one and on the other nodes is zero. Right: histopolant polynomials, the integral of the basis polynomials over the corresponding shaded area evaluates to one and to zero on the others.}
					\label{fig::basis_polynomials}
				\end{center}
			\end{figure}

			It can be shown, \cite{robidoux-polynomial,Gerritsma11}, that if $q_h(\xi)$ is expanded in terms of nodal polynomials, as in
			\eqref{eq::nodal_polynomial_expansion}, then the expansion of its derivative $\frac{\mathrm{d}q_h(\xi)}{\mathrm{d}\xi}$ in terms of histopolant polynomials is
			\begin{align}
				 \left(\frac{\mathrm{d}q_h(\xi)}{\mathrm{d}\xi}\right)_{h} & =
				\sum_{i=1}^{p} \left(\int_{\xi_{i-1}}^{\xi_{i}}\frac{\mathrm{d}q_h(\xi)}{\mathrm{d}\xi}\mathrm{d}\xi\right)e^{p}_{i}(\xi) =
				\sum_{i=1}^{p} \left(q_h(\xi_{i}) - q_h(\xi_{i-1})\right)e^{p}_{i}(\xi) \nonumber \\
 &= \sum_{i=1}^{p} \left(q_{i} - q_{i-1}\right)e^{p}_{i}(\xi) = \sum_{i=1,j=0}^{p}\mathsf{E}^{1,0}_{i,j}q_{j}e^{p}_{i}(\xi)\;,
			\end{align}
			where $\mathsf{E}^{1,0}_{i,j}$ are the coefficients of the $p\times(p+1)$ matrix $\boldsymbol{\mathsf{E}}^{1,0}$
			\begin{equation}
				\boldsymbol{\mathsf{E}}^{1,0} :=
					\left(
						\begin{array}{cccccc}
							-1 & 1 & 0 & 0 &  \dots & 0  \\
							0 & -1 & 1 & 0 & \ddots & 0 \\
							\vdots &    &  \ddots & \ddots  &  & \vdots \\
							0 & \ddots & 0 & -1 & 1 & 0 \\
							0 & \dots & 0 & 0 & -1 & 1
						\end{array}
					\right)\;.
			\end{equation}
			The following identity holds (Commuting property)
			\begin{equation}
				 \left(\frac{\mathrm{d}q(\xi)}{\mathrm{d}\xi}\right)_{h} = \frac{\mathrm{d}q_{h}(\xi)}{\mathrm{d}\xi}\,.
			\end{equation}
			For an example of the one-dimensional basis polynomials corresponding to $p=4$, see \figref{fig::basis_polynomials}.
			
\subsubsection{Two dimensional basis functions}
			Consider the nodal polynomials \eqref{eq::lagrange_interpolants}, $l^{p}_{i}(\xi)$ of degree $p$, the histopolant polynomials \eqref{eq::histopolant_polynomials_definition}, $e^{p}_{i}(\xi)$ of degree $p-1$, the canonical square $\Omega=I\times I\subset \mathbb{R}^{2}$, and take $\xi,\eta\in I=[-1,1]$.
			
			\textbf{Basis functions for $W_{h}$}
				Combining nodal polynomials we can construct the polynomial basis functions for $W_{h}$ on a reference quadrilateral. Consider the canonical interval $I=[-1,1]$, the canonical square $\Omega=I\times I\subset \mathbb{R}^{2}$, the nodal polynomials \eqref{eq::lagrange_interpolants}, $l^{p}_{i}(\xi)$ of degree $p$, and take $\xi,\eta\in I$. Then a set of two-dimensional basis polynomials, $\epsilon^{\,W}_{k}(\xi,\eta; p)$, on $\Omega$ can be constructed as the tensor product of the one-dimensional ones
			\begin{equation}
				\epsilon^{\,W}_{k}(\xi,\eta; p) := l^{p}_{i}(\xi)\,l^{p}_{j}(\eta)\,, \qquad i,j=0,\dots,p, \quad k = j+1 + i(p+1)\,. \label{eq::2d_nodal_polynomials_definition}
			\end{equation}
These polynomials, $\epsilon^{\,W}_{k}(\xi,\eta;p)$, have degree $p$ in each direction and from \eqref{eq::nodal_basis_polynomials} it follows that they satisfy, see
			\cite{Gerritsma11,Palha2014},
			\begin{equation}
				\epsilon^{\,W}_{k}(\xi_{i},\eta_{j}; p) =
				\begin{cases}
					1\, & \mbox{if }k=i(p+1) + j + 1 \\
					 & \\
					0\, & \mbox{if } k\neq i(p+1) + j + 1
				\end{cases}\,, \quad
				i,j = 0,\dots,p, \quad k = 1,\dots,(p+1)^2\,. \label{eq::2d_nodal_polynomials_properties}
			\end{equation}
Where, as before, $\xi_{i}$ and $\eta_{i}$ with $i=0,\dots,p$ are the Gauss-Lobatto-Legendre (GLL) nodes. Let $\omega_h(\xi,\eta)$ be a polynomial function of degree $p$ in $\xi$ and $\eta$, defined on $\Omega$ and
			\begin{equation}
				\omega_{k}^{p} = \omega_h(\xi_{i},\eta_{j}),\qquad \text{with } k = i(p+1) + j+1, \label{eq:2d_nodal_definition_expansion_coefficients}
			\end{equation}
			then its expansion in terms of these polynomials, $\omega_{h}(\xi,\eta;p)$, is given by
			\begin{equation}
				\omega_{h}(\xi,\eta;p) = \sum_{k=1}^{(p+1)^2}\omega_{k}^{p}\,\epsilon_{k}^{\,W}(\xi,\eta;p)\,. \label{eq::2d_nodal_polynomials_expansion}
			\end{equation}
For this relation between the expansion coefficients and nodal interpolation we denote the polynomials in \eqref{eq::2d_nodal_polynomials_definition} by 
\emph{nodal polynomials}. Therefore we set $W_{h}^{p}:= \mathrm{span}\{\epsilon_{1}^{\,W}(\xi,\eta;p),\dots,\epsilon_{(p+1)^{2}}^{\,W}(\xi,\eta;p)\}$. To 
simplify the notation, the explicit reference to the polynomial degree $p$ will be dropped from the function space, the basis functions, and the coefficient 
expansion, therefore from here on we will simply use $W_{h}$, $\epsilon_{k}^{\,W}(\xi,\eta)$, and $\omega_{k}$.
			
			\textbf{Basis functions for $U_{h}$}
In a similar way, but combining nodal polynomials with histopolant polynomials, we can construct the polynomial basis functions for $U_{h}$ on quadrilaterals. Consider the nodal polynomials \eqref{eq::lagrange_interpolants}, $l^{p}_{i}(\xi)$ of degree $p$, the histopolant polynomials \eqref{eq::histopolant_polynomials_definition}, $e^{p}_{i}(\xi)$ of degree $p-1$, the canonical square $\Omega=I\times I\subset \mathbb{R}^{2}$, and take $\xi,\eta\in I=[-1,1]$. A set of two-dimensional basis polynomials, $\vec{\epsilon}^{\,U}_{k}(\xi,\eta;p)$, can be constructed as the tensor product of the one-dimensional basis functions
			\begin{equation}
				\vec{\epsilon}^{\,U}_{k}(\xi,\eta;p):=
				\begin{cases}
					 l_{i}^{p}(\xi)e_{j}^{p}(\eta)\,\vec{e}_{\xi} & \mbox{if } k \leq p(p+1), \quad \mbox{with } i=0,\dots,p, \quad j=1,\dots,p, \quad k = ip + j\;, \\
					& \\
					 e_{i}^{p}(\xi)l_{j}^{p}(\eta)\,\vec{e}_{\eta} & \mbox{if } k > p(p+1), \quad \mbox{with } i=1,\dots,p, \quad j=0,\dots,p, \quad k = (p+i-1)(p+1) + j+1\;.
				\end{cases} \label{eq::edge_polynomials_definition}
			\end{equation}
These polynomials, $\vec{\epsilon}^{\,U}_{k}(\xi,\eta;p)$, have degree $p$ in $\xi$ and $p-1$ in $\eta$ if $k \leq p(p+1)$. If $k> p(p+1)$, then the degree 
in $\xi$ is $p-1$ and the degree in $\eta$ is $p$. Using \eqref{eq::nodal_basis_polynomials} and \eqref{eq::histopolant_polynomials_properties} it follows 
that these polynomials satisfy, \cite{Gerritsma11,Palha2014}
			\begin{equation}
				\int_{\eta_{j-1}}^{\eta_{j}} \vec{\epsilon}_{k}^{\,U}(\xi_{i},\eta;p)\cdot\vec{e}_{\xi}\,\mathrm{d}\xi =
				\begin{cases}
					1 & \mbox{if } k =ip + j \,, \\
					 & \\
					0 & \mbox{if } k \neq ip + j\,,
				\end{cases}\quad\mbox{with}\quad
				\begin{cases}
					i = 0,\dots,p, \\
					j = 1, \dots, p, \\
					k = 1,\dots,2p(p+1)\,,
				\end{cases}
				\label{eq::edge_polynomials_properties_xi}
			\end{equation}

			and
			\begin{equation}
				\int_{\xi_{i-1}}^{\xi_{i}} \vec{\epsilon}_{k}^{\,U}(\xi,\eta_{j})\cdot\vec{e}_{\eta}\,\mathrm{d}\eta =
				\begin{cases}
					1 & \mbox{if } k = (p+i-1)(p+1) + j + 1\,,\\
					 & \\
					0 & \mbox{if } k \neq (p+i-1)(p+1) + j + 1\,,
				\end{cases} \quad\mbox{with}\quad
				\begin{cases}
					i = 1,\dots,p, \\
					j = 0, \dots, p, \\
					k = 1,\dots,2p(p+1)\,.
				\end{cases}
				\label{eq::edge_polynomials_properties_eta}
			\end{equation}
			Let $\vec{u}_h(\xi,\eta;p)$ be a vector valued polynomial function defined on $\Omega$ and
			\begin{equation}
				u_{k}^{p} =
				\begin{dcases}
					\int_{\eta_{j-1}}^{\eta_{j}} \vec{u}_h(\xi_{i},\eta)\cdot\vec{e}_{\xi}\,\mathrm{d}\eta\, & \mbox{if } k \leq p(p+1), \quad \mbox{with} \quad \begin{cases} i=0,\dots,p, \\ j=1,\dots,p, \\ k = ip + j\,, \end{cases} \\
					& \\
					\int_{\xi_{i-1}}^{\xi_{i}} \vec{u}_h(\xi,\eta_{j})\cdot\vec{e}_{\eta}\,\mathrm{d}\xi & \mbox{if } k > p(p+1), \quad \mbox{with} \quad \begin{cases} i=1,\dots,p, \\ j=0,\dots,p, \\ k = (p+i-1)(p+1) + j + 1\,. \end{cases}
				\end{dcases} \label{eq::edge_basis_expansion_coefficientes}
			\end{equation}
			then its expansion in terms of these polynomials, $\vec{u}_{h}(\xi,\eta;p)$, is given by
			\begin{equation}
				\vec{u}_{h}(\xi,\eta;p) = \sum_{k=1}^{2p(p+1)} u_{k}^{p}\vec{\epsilon}^{\,Q}_{k}(\xi,\eta;p)\,. \label{eq::expansion_edge_polynomials}
			\end{equation}
			The expansion $\vec{u}_{h}(\xi,\eta;p)$ is a two-dimensional  polynomial edge histopolant (interpolates integral values along lines).
Since the coefficients of this expansion are edge (or flux) integrals, we denote the polynomials in \eqref{eq::edge_polynomials_definition} by \emph{edge polynomials}.
We set $U_{h} := \mathrm{span}\{\epsilon_{1}^{\,U}(\xi,\eta;p),\dots,\epsilon_{2p(p+1)}^{\,U}(\xi,\eta;p)\}$. To simplify the notation, the explicit reference to the polynomial degree $p$ will be dropped from the function space, the basis functions, and the expansion coefficients, therefore from here on we will simply use $U_{h}$, $\vec{\epsilon}_{k}^{\,U}(\xi,\eta)$, and $u_{k}$.
			
			\textbf{Basis functions for $Q_{h}$}
			Combining histopolant polynomials we can construct the polynomial basis functions for $Q_{h}$ on a quadrilateral. Consider the canonical interval $I=[-1,1]$, the canonical square $\Omega=I\times I\subset \mathbb{R}^{2}$, the histopolant polynomials \eqref{eq::histopolant_polynomials_definition}, $e^{p}_{i}(\xi)$ of degree $p-1$, and take $\xi,\eta\in I$. Then a set of two-dimensional basis polynomials, $\epsilon^{Q}_{k}(\xi,\eta; p)$, can be constructed as the tensor product of the one-dimensional ones
			\begin{equation}
				\epsilon^{Q}_{k}(\xi,\eta; p) := e^{p}_{i}(\xi)\,e^{p}_{j}(\eta), \qquad i,j=1,\dots,p, \quad k = j + (i-1)p\,. \label{eq::volume_polynomials_definition}
			\end{equation}
			These polynomials, $\epsilon^{Q}_{k}(\xi,\eta;p)$, have degree $p-1$ in each variable and satisfy, see \cite{Gerritsma11,Palha2014},
			\begin{equation}
				 \int_{\xi_{i-1}}^{\xi_{i}}\int_{\eta_{j-1}}^{\eta_{j}}\epsilon^{Q}_{k}(\xi,\eta; p)\,\mathrm{d}\xi\mathrm{d}\eta =
				\begin{cases}
					1 & \mbox{if }k=(i-1)p + j \\
					 & \\
					0 & \mbox{if } k\neq (i-1)p + j
				\end{cases}\,, \quad
				i,j = 1,\dots,p, \quad k = 1,\dots,p^{2}\,. \label{eq::volume_polynomials_properties}
			\end{equation}
			Where, as before, $\xi_{i}$ and $\eta_{i}$ with $i=0,\dots,p$ are the Gauss-Lobatto-Legendre (GLL) nodes. Let $K_h(\xi,\eta)$ be a polynomial function defined on $\Omega$ and $K_{k}^{p} = \int_{\xi_{i-1}}^{\xi_{i}}\int_{\eta_{j-1}}^{\eta_{j}}K_h(\xi,\eta)\,\mathrm{d}\xi\mathrm{d}\eta$ with $k = j + (i-1)p$, then its expansion in terms of these polynomials, $K_{h}(\xi,\eta;p)$, is given by
			\begin{equation}
				K_{h}(\xi,\eta;p) = \sum_{k=1}^{p^{2}}K_{k}^{p}\,\epsilon_{k}^{Q}(\xi,\eta;p)\,. \label{eq::volume_polynomials_expansion}
			\end{equation}
For this relation between the expansion coefficients and surface integration we denote the polynomials in \eqref{eq::volume_polynomials_definition} by \emph{surface polynomials}. Moreover, these basis polynomials satisfy $\epsilon^{Q}_{k}(\xi,\eta;p)\in L^{2}(\Omega)$. Therefore we set $Q_{h}^{p}:= \mathrm{span}\{\epsilon_{1}^{Q}(\xi,\eta;p),\dots,\epsilon_{p^{2}}^{Q}(\xi,\eta;p)\}$. To simplify the notation, the explicit reference to the polynomial degree $p$ will be dropped from the function space, the basis functions, and the coefficient expansion, therefore from here on we will simply use $Q_{h}$, $\epsilon_{k}^{Q}(\xi,\eta)$, and $K_{k}$.
			
			\subsubsection{Properties of the basis functions}
The first property that can be shown, \cite{Gerritsma11,Palha2014}, is that if $\omega_h(\xi,\eta) \in W_h$, then
$\nabla^{\perp} \omega_h(\xi,\eta) \in U_h$, where $\omega_h(\xi,\eta)$ is expanded as \eqref{eq::2d_nodal_polynomials_expansion}.

				\begin{align}
				 \left(\nabla^{\perp}\omega_h(\xi,\eta)\right) &\stackrel{\phantom{\eqref{eq:2d_nodal_definition_expansion_coefficients}}}{=} \begin{aligned}[t]
									& \sum_{i=0, j=1}^{p}\left(\int_{\eta_{j-1}}^{\eta_{j}}\nabla^{\perp}\omega_h(\xi_{i},\eta)\cdot\vec{e}_{\xi}
										 \mathrm{d}\eta\right)\vec{\epsilon}_{ip+j}^{\,U}(\xi,\eta) \nonumber \\
									&  \qquad\qquad +  \sum_{i=1,j=0}^{p}\left(\int_{\xi_{i-1}}^{\xi_{i}}\nabla^{\perp}\omega_h(\xi,\eta_{j})\cdot\vec{e}_{\eta}
										 \mathrm{d}\xi\right)\vec{\epsilon}_{(p+i-1)(p+1) + j+1}^{\,U}(\xi,\eta) \nonumber
				   \end{aligned} \nonumber \\
				 & \stackrel{\phantom{\eqref{eq:2d_nodal_definition_expansion_coefficients}}}{=}  \begin{aligned}[t]
				          & \sum_{i=0,j=1}^{p}\left(\omega_h(\xi_{i},\eta_{j-1})-\omega_h(\xi_{i},\eta_{j})\right)\vec{\epsilon}_{ip+j}^{\,U}(\xi,\eta) \nonumber \\
				          & \qquad\qquad +  \sum_{i=1,j=0}^{p}\left(\omega_h(\xi_{i},\eta_{j}) - \omega_h(\xi_{i-1},\eta_{j})\right)
							 \vec{\epsilon}_{(p+i-1)(p+1)+j+1}^{\,U}(\xi,\eta)\nonumber\\
				          \end{aligned}\nonumber \\
				 & \stackrel{\eqref{eq:2d_nodal_definition_expansion_coefficients}}{=}  \sum_{i=0,j=1}^{p}\left(\omega_{ip+j}-\omega_{ip+j+1}\right)
						\vec{\epsilon}_{ip+j}^{\,U}(\xi,\eta) +  \sum_{i=1,j=0}^{p}\left(\omega_{i(p+1) + j+1} - \omega_{(i-1)(p+1)+j+1}\right)
						 \vec{\epsilon}_{(p+i-1)(p+1)+j+1}^{\,U}(\xi,\eta)\nonumber\\
				 &\stackrel{\phantom{\eqref{eq:2d_nodal_definition_expansion_coefficients}}}{=} 
\sum_{k=1}^{2p(p+1)} \sum_{j=1}^{(p+1)^2}
\mathsf{E}^{1,0}_{k,j} \omega_{j} \vec{\epsilon}_{k}^{\,U}(\xi,\eta)\,, \label{eq::2d_nodal_basis_expansion_curl}
			\end{align}
			where $\mathsf{E}_{k,j}^{1,0}$ are the coefficients of the $2p(p+1)\times (p+1)^{2}$ matrix $\boldsymbol{\mathsf{E}}^{1,0}$ and are defined as
			\begin{equation}
				\mathsf{E}_{k,j}^{1,0} :=
				\begin{cases}
					-1 & \mbox{if } j = k + 1 + (k\,\mathrm{div}(p+1))\quad \text{and} \quad 1 \leq k \leq p(p+1)\,, \\
					1  & \mbox{if } j = k + (k\,\mathrm{div}(p+1))\quad \text{and} \quad 0 \leq k \leq p(p+1)\,, \\
					1  & \mbox{if } j = k - (p-1)(p+1)\quad \text{and} \quad p(p+1) < k \leq 2p(p+1)\,,\\
					-1 & \mbox{if } j = k - p(p+1)\quad \text{and} \quad p(p+1) < k \leq 2p(p+1)\,, \\
					0  & \mbox{otherwise}  \,.
				\end{cases}
			\end{equation}
			Here $(k\,\mathrm{div} p)$ denotes integer division in which the remainder is discarded.

From \eqref{eq::2d_nodal_basis_expansion_curl} it follows that
\[ \omega_{h}(\xi,\eta) \in W_h \quad \Longrightarrow \quad \nabla^{\perp} \omega_h(\xi,\eta) \in U_h \;,\]
which is the finite-dimensional analogue of
\[ \omega(\xi,\eta) \in H(\mathrm{rot};\Omega) \quad \Longrightarrow \quad \nabla^{\perp} \omega_(\xi,\eta) \in H(\mathrm{div};\Omega) \;.\]
Or fully discrete, if $\omega_j$ are the expansion coefficients of $\omega_h \in W_h$ with respect to the basis $e_i^W(\xi,\eta)$, then
$\mathsf{E}^{1,0}_{k,j} \omega_{j}$ are the expansion coefficients of $\nabla^{\perp} \omega_h$ in $U_h$ with respect to the basis $\vec{u}_j^U(\xi,\eta)$.
			
			As a special case we have that
			\begin{equation}
				\nabla^{\perp} \epsilon^{\,W}_{j} = \sum_{k=1}^{2p(p+1)}\mathsf{E}_{k,j}^{1,0}\vec{\epsilon}^{\,U}_{k}\,, \label{eq:hilbert_subcomplex_basis_W}
			\end{equation}
			and therefore $\nabla^{\perp}\vec{\epsilon}^{\,W}_{j}\in U_{h}$, with $j=1,\dots,(p+1)^{2}$, and these basis functions satisfy \eqref{eq:hilbert_subcomplex}.
			
				The second property that can be shown, \cite{Gerritsma11,Palha2014}, is that if $\vec{u}_h(\xi,\eta)\in U_h$ is expanded in terms of edge polynomials, as in \eqref{eq::expansion_edge_polynomials}, then the expansion of $\nabla\cdot\vec{u}_h$ in terms of the surface polynomials, \eqref{eq::volume_polynomials_definition}, is
			\begin{align}
				\nabla\cdot\vec{u}_h(\xi,\eta) &\stackrel{\phantom{\eqref{eq::edge_basis_expansion_coefficientes}}}{=} \sum_{i,j=1}^{p}\left(\int_{\xi_{i-1}}^{\xi_{i}}\int_{\eta_{j-1}}^{\eta_{j}}\nabla\cdot\vec{u}_h(\xi,\eta)\,\mathrm{d}\xi\mathrm{d}\eta\right)\epsilon_{j+(i-1)p}^{Q}(\xi,\eta) \nonumber\\
				 &\stackrel{\phantom{\eqref{eq::edge_basis_expansion_coefficientes}}}{=}\begin{aligned}[t]
				 		 &\sum_{i,j=1}^{p}\left(\int_{\eta_{j-1}}^{\eta_{i}}\vec{u}_h(\xi_{i},\eta)\cdot\vec{e}_{\xi}\,\mathrm{d}\eta + \int_{\xi_{i-1}}^{\xi_{i}}\vec{u}_h(\xi,\eta_{j})\cdot\vec{e}_{\eta}\,\mathrm{d}\xi -\int_{\eta_{j-1}}^{\eta_{i}}\vec{u}_h(\xi_{i-1},\eta)\cdot\vec{e}_{\xi}\,\mathrm{d}\eta\right. \nonumber \\
						&\qquad\qquad\left.- \int_{\xi_{i-1}}^{\xi_{i}}\vec{u}_h(\xi,\eta_{j-1})\cdot\vec{e}_{\eta}\,\mathrm{d}\xi\right)\epsilon_{j+(i-1)p}^{Q}(\xi,\eta)
					\end{aligned} \nonumber\\
				 &\stackrel{\eqref{eq::edge_basis_expansion_coefficientes}}{=} \sum_{k=1}^{p^{2}} \left(u_{k+p} + u_{k+p(p+1)+(k\mathrm{div} (p+1))+1} - u_{k} - u_{k+p(p+1) + (k\mathrm{div}(p+1))} \right)\epsilon_{k}^{Q}(\xi,\eta)\nonumber \\
				 &\stackrel{\phantom{\eqref{eq::edge_basis_expansion_coefficientes}}}{=} \sum_{k=1}^{p^{2}} \sum_{j=1}^{2p(p+1)}\mathsf{E}^{2,1}_{k,j} u_{j} \epsilon_{k}^{Q}(\xi,\eta)\,, \label{eq::edge_basis_expansion_curl}
			\end{align}
			where $\mathsf{E}^{2,1}_{k,j}$ are the coefficients of the $p^{2}\times 2p(p+1)$ incidence matrix $\boldsymbol{\mathsf{E}}^{2,1}$ defined as
			\begin{equation}
				\mathsf{E}_{k,j}^{2,1} :=
				\begin{cases}
					1  & \mbox{if } j = k + p\,, \\
					1  & \mbox{if } j = k + p(p+1) + (k\,\mathrm{div}(p+1))+1\,, \\
					-1 & \mbox{if } j = k \,, \\
					-1 & \mbox{if } j = k + p(p+1) + (k\,\mathrm{div}(p+1))\,, \\
					0  & \mbox{otherwise}\,.
				\end{cases}
			\end{equation}
Equation \eqref{eq::edge_basis_expansion_curl} confirms that we have a finite dimensional Hilbert sequence as in \eqref{eq:hilbert_subcomplex}, because
\[ \vec{u}_h(\xi,\eta)  \in U_h \quad \Longrightarrow \quad \nabla \cdot \vec{u}_h(\xi,\eta) \in Q_h \;,\]
which is the finite dimensional analogue of
\[ \vec{u} \in H(\mathrm{div};\Omega) \quad \Longrightarrow \quad \nabla \cdot \vec{u} \in L^2(\Omega)\;.\]
In terms of the expansion coefficients we have: If $u_j$ are the expansion coefficients of $\vec{u}_h \in U_h$ with respect to the basis $\vec{\epsilon}_j^U$, then the expansion coefficients of $\nabla \cdot \vec{u}_h \in Q_h$ with respect to the basis $\epsilon_k^Q$ are given by $\sum_{j=1}^{2p(p+1)}\mathsf{E}^{2,1}_{k,j} u_{j}$.
			As a special case we have that
			\begin{equation}
				\nabla\cdot\vec{\epsilon}^{\,U}_{j} = \sum_{k=1}^{p^{2}}\mathsf{E}_{k,j}^{2,1}\epsilon^{Q}_{k}\,. \label{eq:hilbert_subcomplex_basis_U}
			\end{equation}
If $\omega_j$ are the expansion coefficients of $\omega_h \in W_h$, then $\mathsf{E}_{k,j}^{1,0} \omega_j$ are the expansion coefficients of
$\nabla^{\perp} \omega_h \in U_h$. Then $\mathsf{E}_{i,k}^{2,1}\mathsf{E}_{k,j}^{1,0} \omega_j$ are the expansion coefficients of
$\nabla \cdot \nabla^{\perp} \omega_h \in Q_h$. Since $\nabla \cdot \nabla^{\perp} \omega_h =0$ for all $\omega_h$ (and therefore
$\omega_j$) and because $\epsilon_k^Q$ forms a basis for $Q_h$, we need to have that \[ \mathsf{E}_{i,k}^{2,1} \circ \mathsf{E}_{k,j}^{1,0} \equiv 0 \;.\]
This is the fully discrete representation of the vector identity $\nabla \cdot \nabla^{\perp} \equiv 0$.
			
		\subsubsection{Some remarks}	
Note that the degrees of freedom for $W_h$ are associated with the points in the GLL-grid. In a multi-element setting, neighboring elements share
the GLL-points on the boundary of the element, thus imposing $C^0$-continuity for functions in $W_h$.
The degrees of freedom for $U_h$ are the integrals of the normal components (the flux) of $\vec u_h$ over the edge in the GLL-grid. In a multi-element
setting, neighboring elements share an edge and therefore also the degree of freedom associated with that edge. So only the normal component of
$\vec{u}_h \in U_h$ is continuous between elements, thus making $U_h$ a proper subspace of $H(\mathrm{div};\Omega)$.
The degrees of freedom for $Q_h$ are the integral values over the two-dimensional surfaces in the GLL-grid. A surface in one spectral element is
entirely disjoint from a surface in neighboring elements and therefore, the approximation in $Q_h$ is discontinuous between elements. In that sense
$Q_h$ is a proper subspace of $L^2(\Omega)$.

There is a close relation between this mimetic finite element discretization and the traditional C-grid finite difference discretization \cite{AL81},
where rotational moments are located at vertices, normal velocities on edges
and pressure and mass variables at cell centers. It is possible to construct an analogue to the D-grid finite difference discretization, for which
the tangential and not the normal velocities reside on the edges. Two advantages of the mimetic finite element discretization is the immediate
generalization to arbitrary order and the possibility to treat deformed meshes, which we have not addressed in the current work.

\firstRev{There exists for the Hilbert subcomplex of discrete function spaces described above a discrete Helmholtz decomposition \cite{CS12} 
of vector fields in $U_h$ into a unique partition of rotational components in $W_h$ and divergent components in $Q_h$ which are orthogonal 
to one another. This result has been shown previously for the mixed mimetic spectral element function spaces \cite{KPG11}.}

	\subsection{Discrete weak formulation}
		Consider again the domain $\Omega\subset\mathbb{R}^{2}$  and its tessellation $\mathcal{T}(\Omega)$ consisting of $N$ arbitrary quadrilaterals (possibly curved), $\Omega_{m}$, with $m = 1, \dots, N$. We assume that all quadrilateral elements $\Omega_{m}$ can be obtained from a map $\Phi_{m}:(\xi,\eta)\in I^{2} \mapsto (x,y)\in\Omega_{m}$. Then the pushforward $\Phi_{m,*}$ maps functions in the reference element $I^{2}$ to functions in the physical element $\Omega_{m}$, see for example \cite{abraham_diff_geom,frankel}. For this reason it suffices to explore the analysis on the reference domain $I^{2}$. Additionally, the multi-element case follows the standard approach in finite elements.
		
\begin{remark}
	If a differential geometry formulation was used, the physical quantities would be represented by differential $k$-forms and the map $\Phi_{m}:(\xi,\eta)\in I^{2} \mapsto(x,y)\in\Omega_{m}$ would generate a pullback, $\Phi_{m}^{*}$, mapping $k$-forms in physical space, $\Omega_{m}$, to $k$-forms in the reference element, $I^{2}$, \cite{KPG11}.
\end{remark}
		
		The discrete weak formulation can be stated as: given $\Omega = I^{2}$, the polynomial degree $p$ and a Coriolis term $f_{h}\in W_{h}$, for any time $t\in (0,t_{F}]$ find $\vec{u}_{h},\vec{F}_{h}\in U_{h}$, $h_{h},K_{h}\in Q_{h}$, and $q_{h}\in W_{h}$ such that
	 \begin{subnumcases}{\label{eq:shallow_water_continuous_weak_form_discrete}}
		\ip{\frac{\partial\vec{u}_{h}}{\partial t}}{\vec{\nu}_{h}} + \ip{q_{h}\times\vec{F}_{h}}{\vec{\nu}_{h}} - \ip{K_{h} + gh_{h}}{\nabla\cdot\vec{\nu}_{h}} = 0\,, &
		$\forall \vec{\nu}_{h}\in U_{_{h}}$  \label{mom_cont_weak_discrete} \\
		\ip{\frac{\partial h_{h}}{\partial t}}{\sigma_{h}} + \ip{\nabla\cdot\vec{F}_{h}}{\sigma_{h}} = 0\,, & $\forall \sigma_{h}\in Q_{h}$, \label{mas_cont_weak_discrete} \\
		\ip{h_{h}q_{h}}{\zeta_{h}} = -\ip{\vec u_{h}}{\nabla^{\perp}\zeta_{h}} + \ip{f_{h}}{\zeta_{h}}\,, & $\forall \zeta_{h}\in W_{h}$, \label{eq:definition_potential_vorticity_weak_discrete} \\
		\ip{\vec{F}_{h}}{\vec{\varphi}_{h}} = \ip{h_{h}\vec{u}_{h}}{\vec{\varphi}_{h}}\,, &  $\forall \vec{\varphi}_{h}\in U_{h}$, \label{eq:definition_h_flux_weak_discrete} \\
		\ip{K_{h}}{\kappa_{h}} = \frac{1}{2}\ip{\vec{u}_{h}\cdot\vec{u}_{h}}{\kappa_{h}}\,, & $\forall \kappa_{h}\in Q_{h}$.  \label{eq:definition_kinetic_energy_weak_discrete}
	\end{subnumcases}
	
	Using the expansions for $\vec{u}_{h}$, $\vec{F}_{h}$, $h_{h}$, $K_{h}$ and $q_{h}$ in \eqref{eq:physical_quantities_expansion},
	\eqref{eq:shallow_water_continuous_weak_form_discrete} can be written as: find $\boldsymbol{u},\boldsymbol{F}\in\mathbb{R}^{d_{U}}$,
	$\boldsymbol{h},\boldsymbol{K}\in\mathbb{R}^{d_{Q}}$, and $\boldsymbol{q}\in\mathbb{R}^{d_{W}}$ such that
	 \begin{subnumcases}{\label{eq:shallow_water_continuous_weak_form_discrete_expansion}}
		 \sum_{i=1}^{d_{U}}\ip{\vec{\epsilon}^{\,U}_{i}}{\vec{\epsilon}^{\,U}_{j}}\frac{\mathrm{d}u_{i}}{\mathrm{d}t} +
		 \sum_{i=1}^{d_{U}}\ip{q_{h}\times\vec{\epsilon}_{i}^{\,U}}{\vec{\epsilon}_{j}^{\,U}}F_i -
		 \sum_{i,k=1}^{d_{Q}}\mathsf{E}^{2,1}_{k,j}\ip{\epsilon_{i}^{Q}}{\epsilon_{k}^{Q}}(K_i + gh_i) = 0\,, & $j = 1,\dots, d_{U}\,,$  \label{mom_cont_weak_discrete_expansion} \\
		 \sum_{i=1}^{d_{Q}}\ip{\epsilon_{i}^{Q}}{\epsilon_{j}^{Q}}\frac{\mathrm{d}h_{i}}{\mathrm{d}t} +
		\sum_{i=1}^{d_{U}} \sum_{i,k=1}^{d_{Q}}\mathsf{E}^{2,1}_{k,i}\ip{\epsilon_{k}^{Q}}{\epsilon_{j}^{Q}}F_i = 0\,, & $j=1,\dots,d_{Q}$\,, \label{mas_cont_weak_discrete_expansion} \\
		 \sum_{i=1}^{d_{W}}\ip{h_{h}\epsilon_{i}^{\,W}}{\epsilon_{j}^{\,W}}q_i = -\sum_{i,k=1}^{d_{U}} \mathsf{E}_{k,j}^{1,0}\ip{\vec{\epsilon}_{i}^{\,U}}{\vec{\epsilon}_{k}^{\,U}}u_i +
		 \sum_{i=1}^{d_{W}}\ip{\epsilon_{i}^{\,W}}{\epsilon_{j}^{\,W}}f_i\,, & $j=1,\dots,d_{W}$\,, \label{eq:definition_potential_vorticity_weak_discrete_expansion} \\
		 \sum_{i=1}^{d_{U}}\ip{\vec{\epsilon}_{i}^{\,U}}{\vec{\epsilon}_{j}^{\,U}}F_i =
		 \sum_{i=1}^{d_U}\ip{h_{h}\vec{\epsilon}_{i}^{\,U}}{\vec{\epsilon}_{j}^{\,U}}u_i\,, &  $j=1,\dots,d_{U}$\,, \label{eq:definition_h_flux_weak_discrete_expansion} \\
		 \sum_{i=1}^{d_{Q}}\ip{\epsilon_{i}^{Q}}{\epsilon_{j}^{Q}}K_i =
		 \frac{1}{2}\sum_{i=1}^{d_{U}}\ip{\vec{u}_{h}\cdot\vec{\epsilon}_{i}^{\,U}}{\epsilon_{j}^{Q}}u_i\,, & $j=1,\dots,d_{Q}$\,.  \label{eq:definition_kinetic_energy_weak_discrete_expansion}
	\end{subnumcases}
	with $\boldsymbol{u}:=[u_{1},\dots,u_{d_{U}}]^{\top}$, $\boldsymbol{F}:=[F_{1},\dots,F_{d_{U}}]^{\top}$, $\boldsymbol{h}:=[h_{1},\dots,h_{d_{Q}}]^{\top}$,
	$\boldsymbol{K}:=[K_{1},\dots,K_{d_{Q}}]^{\top}$ and $\boldsymbol{q}:=[q_{1},\dots,q_{d_{W}}]^{\top}$.

	Using matrix notation, \eqref{eq:shallow_water_continuous_weak_form_discrete_expansion} can be expressed more compactly as
	 \begin{subnumcases}{\label{eq:shallow_water_continuous_weak_form_discrete_matrix}}
		 \boldsymbol{\mathsf{U}}\frac{\mathrm{d}\boldsymbol{u}}{\mathrm{d}t} + \boldsymbol{\mathsf{U}}^{q}\boldsymbol{F} - \left(\boldsymbol{\mathsf{E}}^{2,1}\right)^{\top}\boldsymbol{\mathsf{Q}}\,(\boldsymbol{K}+g\boldsymbol{h}) = 0\,, &  \label{mom_cont_weak_discrete_matrix} \\
		 \boldsymbol{\mathsf{Q}}\, \frac{\mathrm{d}\boldsymbol{h}}{\mathrm{d}t} + \boldsymbol{\mathsf{Q}}\, \boldsymbol{\mathsf{E}}^{2,1}\boldsymbol{F} = 0\,, & \label{mas_cont_weak_discrete_matrix} \\
		\boldsymbol{\mathsf{W}}^{h}\boldsymbol{q} = -\left(\boldsymbol{\mathsf{E}}^{1,0}\right)^{\top}\boldsymbol{\mathsf{U}}\,\boldsymbol{u} + \boldsymbol{\mathsf{W}}\boldsymbol{f}\,, & \label{eq:definition_potential_vorticity_weak_discrete_matrix} \\
		\boldsymbol{\mathsf{U}}\boldsymbol{F} = \boldsymbol{\mathsf{U}}^{h}\boldsymbol{u}\,, &  \label{eq:definition_h_flux_weak_discrete_matrix} \\
		\boldsymbol{\mathsf{Q}}\boldsymbol{K} = \frac{1}{2}\boldsymbol{\mathsf{U}}^{u}\boldsymbol{u}\,. &  \label{eq:definition_kinetic_energy_weak_discrete_matrix}
	\end{subnumcases}
	The coefficients of the matrices $\boldsymbol{\mathsf{U}}$, $\boldsymbol{\mathsf{Q}}$, and $\boldsymbol{\mathsf{W}}$ are given by
	\begin{equation}
		\mathsf{U}_{ij} := \ip{\vec{\epsilon}_{i}^{\,U}}{\vec{\epsilon}_{j}^{\,U}},\quad \mathsf{Q}_{ij} := \ip{\epsilon_{i}^{Q}}{\epsilon_{j}^{Q}}, \quad\mathrm{and}\quad \mathsf{W}_{ij} := \ip{\epsilon_{i}^{\,W}}{\epsilon_{j}^{\,W}}\,.
	\end{equation}
	Similarly, the coefficients of the matrices $\boldsymbol{\mathsf{U}}^{q}$, $\boldsymbol{\mathsf{U}}^{h}$, $\boldsymbol{\mathsf{U}}^{u}$, and $\boldsymbol{\mathsf{W}}^{h}$ are given by
	 \begin{equation}
		\mathsf{U}^{q}_{ij} := \ip{q_{h}\times\vec{\epsilon}_{j}^{\,U}}{\vec{\epsilon}_{i}^{\,U}},\quad \mathsf{U}^{h}_{ij} := \ip{h_{h}\vec{\epsilon}_{j}^{\,U}}{\vec{\epsilon}_{i}^{\,U}}, \quad \mathsf{U}^{u}_{ij} := \ip{\vec{u}_{h}\cdot\vec{\epsilon}_{j}^{\,U}}{{\epsilon}_{i}^{\,Q}}, \quad\mathrm{and}\quad\mathsf{W}^{h}_{ij} := \ip{h_{h}\epsilon_{j}^{\,W}}{\epsilon_{i}^{\,W}}\,.
	\end{equation}

The mass matrix $\boldsymbol{\mathsf{Q}}$ may be canceled entirely from \eqref{mas_cont_weak_discrete_matrix}, such that the continuity equation holds in the
strong form. 
\[ \frac{\mathrm{d}\boldsymbol{h}}{\mathrm{d}t} + \boldsymbol{\mathsf{E}}^{2,1}\boldsymbol{F} = 0\,. \]
This reflects the fact that the divergence theorem is satisfied point wise as given by \eqref{eq:hilbert_subcomplex_basis_U}.
\secondRev{We further note that \eqref{eq:definition_kinetic_energy_weak_discrete}, \eqref{eq:definition_kinetic_energy_weak_discrete_expansion} and 
\eqref{eq:definition_kinetic_energy_weak_discrete_matrix} are redundant since $\vec u_h\cdot\vec u_h$ may be directly projected onto $\vec\nu_h$ in
\eqref{mom_cont_weak_discrete}, \eqref{mom_cont_weak_discrete_expansion} and \eqref{mom_cont_weak_discrete_matrix} respectively. We have preserved this
expression however as it makes explicit that $K_h\in Q_h$ is a scalar quantity in the same discrete function space as $h_h$, and is efficient to compute 
since it may be done so individually for each element since $Q_h$ is discontinuous across element boundaries.}
	
\section{Discrete conservation properties}

	\subsection{Conservation of mass}

\firstRev{
        Conservation of mass is given by
        \begin{equation}
            \ddt{}\int_{\Omega}h\,\mathrm{d}\Omega = 0\,. \label{eq:conservation_mass_definition}
        \end{equation}
        At the discrete level we have that $h_{h}\in Q_{h}$, therefore $h_{h} = \sum_{j}h_{j}\epsilon^{Q}_{j}$. Recalling that the basis functions $\epsilon^{Q}_{j}$ satisfy (40), we have that
        \begin{equation}
            \int_{\Omega}h_{h}\,\mathrm{d}\Omega = \sum_{j}h_{j} = \valgebra{1}^{\transpose}\valgebra{h}\,, \label{eq:total_mass_discrete}
        \end{equation}
        with $\valgebra{1}^{\transpose}=[1,\dots,1]\in\mathbb{R}^{d_{Q}}$. Eliminating $\moperator{Q}$ in (50b) results in
        \begin{equation}
            \ddt{\valgebra{h}} = -\moperator{E}^{2,1}\valgebra{F}\,. \label{eq:time_variation_of_mass_and_flux}
        \end{equation}
        Combining  \eqref{eq:conservation_mass_definition} with \eqref{eq:total_mass_discrete} and \eqref{eq:time_variation_of_mass_and_flux}, results in conservation of mass at the algebraic level
        \begin{equation}
            \ddt{}\int_{\Omega}h_{h}\,\mathrm{d}\Omega \stackrel{\eqref{eq:total_mass_discrete}}{=} \ddt{}\left(\valgebra{1}^{\transpose}\valgebra{h}\right) = \valgebra{1}^{\transpose}\ddt{\valgebra{h}} \stackrel{\eqref{eq:time_variation_of_mass_and_flux}}{=} - \valgebra{1}^{\transpose}\moperator{E}^{2,1}\valgebra{F} = 0\,.
        \end{equation}
        The last identity results from the telescoping property of the incidence matrix on periodic domains, $\valgebra{1}^{\transpose}\moperator{E}^{2,1} = 0$.
}

\subsection{Conservation of vorticity}

\firstRev{
        Conservation of total vorticity is given by
        \begin{equation}
            \ddt{}\int_{\Omega}\omega\,\mathrm{d}\Omega = 0\,. \label{eq:conservation_vorticity_definition}
        \end{equation}
        At the discrete level, vorticity is defined as: for any time $t\in(0,t_{F}]$, given $\vfield{u}_{h}\in U_{h}$, find $\omega_{h}\in W_{h}$ such that
        \begin{equation}
            \innerproduct{\varsigma_{h}}{\omega_{h}} = - \innerproduct{\gradientperp\varsigma_{h}}{\vfield{u}_{h}}\,,\qquad \forall \varsigma_{h}\in W_{h}\,. \label{eq:definition_vorticity_inner_product}
        \end{equation}
        At the algebraic level \eqref{eq:definition_vorticity_inner_product} becomes
        \begin{equation}
            \moperator{W}\valgebra{\omega} = - \left(\moperator{E}^{1,0}\right)^{\transpose}\moperator{U}\valgebra{u}\,,
        \end{equation}
        and its time derivative
        \begin{equation}
            \moperator{W}\ddt{\valgebra{\omega}} = - \left(\moperator{E}^{1,0}\right)^{\transpose}\moperator{U}\ddt{\valgebra{u}}\,. \label{eq:time_derivative_definition_vorticity}
        \end{equation}
        Since $\omega_{h}\in W_{h}$, we have that $\omega_{h} = \sum_{j}\omega_{k}\epsilon_{j}^{W}$. Recalling that $\sum_{j}\epsilon_{j}^{W} = 1$, we have
        \begin{equation}
            \int_{\Omega}\omega_{h}\,\mathrm{d}\Omega = \int_{\Omega}\sum_{j}\omega_{j}\epsilon_{j}^{W}\,\mathrm{d}\Omega = \int_{\Omega}\sum_{j} \left(\sum_{k}\epsilon_{k}^{W}\right)\omega_{j}\epsilon_{j}^{W}\,\mathrm{d}\Omega = \sum_{j,k}\omega_{j}\int_{\Omega}\epsilon_{k}^{W}\epsilon_{j}^{W}\,\mathrm{d}\Omega = \valgebra{1}^{\transpose}\moperator{W}\valgebra{\omega}\,, \label{eq:discrete_total_vorticity}
        \end{equation}
        with $\valgebra{1}^{\transpose} = [1,\dots,1]\in\mathbb{R}^{d_{W}}$. Therefore \eqref{eq:conservation_vorticity_definition} becomes
        \begin{equation}
            \ddt{}\int_{\Omega}\omega_{h}\,\mathrm{d}\Omega \stackrel{\eqref{eq:discrete_total_vorticity}}{=} \ddt{}\left(\valgebra{1}^{\transpose}\moperator{W}\valgebra{\omega}\right) = \valgebra{1}^{\transpose}\moperator{W} \ddt{\valgebra{\omega}}\,, \label{eq:time_variation_vorticity_algebraic}
        \end{equation}
        Multiplying both sides of \eqref{eq:time_derivative_definition_vorticity} by $\valgebra{1}^{\transpose}$ and combining with \eqref{eq:time_variation_vorticity_algebraic} gives the time conservation of total vorticity at the algebraic level
        \begin{equation}
            \ddt{}\int_{\Omega}\omega_{h}\,\mathrm{d}\Omega = \valgebra{1}^{\transpose}\moperator{W} \ddt{\valgebra{\omega}} = - \valgebra{1}^{\transpose}\left(\moperator{E}^{1,0}\right)^{\transpose}\moperator{U}\ddt{\valgebra{u}} = 0\,.
        \end{equation}
        The last identity again results from the telescoping property of the incidence matrix on periodic domains, $\valgebra{1}^{\transpose}\left(\moperator{E}^{1,0}\right)^{\transpose} = \left(\moperator{E}^{1,0}\valgebra{1}\right)^{\transpose} = 0$.
}

Note that the conservation of vorticity is satisfied irrespective of how $q_h$ is constructed in
\eqref{mom_cont_weak_discrete_expansion}. As such vorticity
conservation is preserved in the event that the anticipated potential vorticity method \cite{SB85} is
used to truncate the potential enstrophy cascade, since this involves removing some downstream \emph{anticipated}
potential vorticity from the right hand side in order to introduce some dispersion into the
vorticity advection equation. We construct this anticipated potential vorticity $\hat{q}_h$
in the weak form as

\begin{equation}\label{apv}
\ip{\hat{q}_h}{\varsigma_h} = \ip{q_h}{\varsigma_h} - \Delta\tau\ip{\vec u_h\times\nabla^{\perp}q_h}{\varsigma_h}
\end{equation}
where $\Delta\tau$ is some time scale associated with the evaluation of the downstream potential vorticity.

\subsection{Conservation of energy}

\firstRev{
        For physical problems governed by the shallow water equations, the total energy $\mathcal{E}$, is given by the sum of kinetic and potential energies
        \begin{equation}
            \mathcal{E} := \int_{\Omega}\left(Kh + \frac{1}{2}gh^{2}\right)\,\mathrm{d}\Omega = \innerproduct{h}{K} + \frac{1}{2}\innerproduct{h}{gh}\,.
        \end{equation}
        Conservation of total energy is then expressed as
        \begin{equation}
            \ddt{\mathcal{E}} = \innerproduct{h}{\pddt{K}} + \innerproduct{\pddt{h}}{K} + \innerproduct{\pddt{h}}{gh}\,.
        \end{equation}
        For the discrete variables, the time variation of energy takes a similar form
        \begin{equation}
            \ddt{\mathcal{E}_{h}} = \innerproduct{h_{h}}{\pddt{K_{h}}} + \innerproduct{\pddt{h_{h}}}{K_{h}} + \innerproduct{\pddt{h_{h}}}{gh_{h}}\,,
        \end{equation}
        which can be written at the algebraic level as
        \begin{equation}
            \ddt{\mathcal{E}_{h}} = \valgebra{h}^{\transpose}\moperator{Q}\ddt{\valgebra{K}} + \ddt{\valgebra{h}^{\transpose}}\moperator{Q}\valgebra{K} + g\ddt{\valgebra{h}^{\transpose}}\moperator{Q}\valgebra{h}\,. \label{eq:energy_conservation_algebraic_1}
        \end{equation}
        From (50e) we have that
        \begin{equation}
            \moperator{Q}\valgebra{K} = \frac{1}{2}\moperator{U}^{u}\valgebra{u}\,,
        \end{equation}
        and therefore its time derivative is given by
        \begin{equation}
            \moperator{Q}\ddt{\valgebra{K}} = \frac{1}{2}\ddt{}\left(\moperator{U}^{u}\valgebra{u}\right)\,. \label{eq:time_derivative_kinetic_energy_algebraic}
        \end{equation}

        If we now recall the definition of $\moperator{U}^{u}$, (52), we can rewrite the right hand side of \eqref{eq:time_derivative_kinetic_energy_algebraic} as
        \begin{equation}
            \frac{1}{2}\ddt{}\left(\moperator{U}^{u}\valgebra{u}\right) = \frac{1}{2}\ddt{}\left(\sum_{j=1}^{d_{U}}\moperatorcomponents{U}^{u}_{ij}\valgebracomponents{u_{j}}\right) = \frac{1}{2}\ddt{}\left(\sum_{k,j=1}^{d_{U}}\valgebracomponents{u}_{k}\innerproduct{\epsilon^{Q}_{i}}{\vfield{\epsilon}_{k}^{\,U}\cdot\vfield{\epsilon}_{j}^{\,U}}\valgebracomponents{u_{j}}\right) = \sum_{j=1}^{d_{U}}\moperatorcomponents{U}^{u}_{ij}\ddt{\valgebracomponents{u_{j}}} = \moperator{U}^{u}\ddt{\valgebra{u}}\,. \label{eq:time_derivative_kinetic_energy_algebraic_2}
        \end{equation}
        Substituting \eqref{eq:time_derivative_kinetic_energy_algebraic_2} into \eqref{eq:energy_conservation_algebraic_1} yields
        \begin{equation}
            \ddt{\mathcal{E}_{h}} = \valgebra{h}^{\transpose}\moperator{U}^{u}\ddt{\valgebra{u}} + \ddt{\valgebra{h}^{\transpose}}\moperator{Q}\valgebra{K} + g\ddt{\valgebra{h}^{\transpose}}\moperator{Q}\valgebra{h}\,. \label{eq:time_derivative_kinetic_energy_algebraic_3}
        \end{equation}
        Noting that
        \begin{equation}
             \valgebra{h}^{\transpose}\moperator{U}^{u} = \sum_{i=1}^{d_{Q}}\valgebracomponents{h}_{i}\moperatorcomponents{U}^{u}_{ij}= \sum_{i,k=1}^{d_{Q},d_{U}}\valgebracomponents{h}_{i}\innerproduct{\epsilon_{i}^{Q}}{\vfield{\epsilon}_{k}^{\,U}\cdot\vfield{\epsilon}_{j}^{\,U}}\valgebracomponents{u}_{k}  = \sum_{i,k=1}^{d_{Q},d_{U}}\valgebracomponents{h}_{i}\innerproduct{\vfield{\epsilon}_{k}^{\,U}}{\epsilon_{i}^{Q}\vfield{\epsilon}_{j}^{\,U}}\valgebracomponents{u}_{k} = \sum_{k=1}^{d_{U}}\valgebracomponents{u}_{k}\moperatorcomponents{U}^{h}_{kj} = \valgebra{u}^{\transpose}\moperator{U}^{h}\,,
        \end{equation}
        we can rewrite \eqref{eq:time_derivative_kinetic_energy_algebraic_3} as
        \begin{equation}
            \ddt{\mathcal{E}_{h}} = \valgebra{u}^{\transpose}\left(\moperator{U}^{h}\right)^{\transpose}\ddt{\valgebra{u}} + \ddt{\valgebra{h}^{\transpose}}\moperator{Q}\valgebra{K} + g\ddt{\valgebra{h}^{\transpose}}\moperator{Q}\valgebra{h}\,,\label{eq:time_derivative_kinetic_energy_algebraic_4}
        \end{equation}
        since $\moperator{U}^{h}$ is a symmetric matrix. If we now use (50d) on the first term on the right hand side of \eqref{eq:time_derivative_kinetic_energy_algebraic_4} and (50b) on the other two terms, we obtain
        \begin{equation}
            \ddt{\mathcal{E}_{h}} = \valgebra{F}^{\transpose}\moperator{U}\ddt{\valgebra{u}} - \valgebra{F}^{\transpose}\left(\moperator{E}^{2,1}\right)^{\transpose}\moperator{Q}(\valgebra{K} + g\valgebra{h})\,,\label{eq:time_derivative_kinetic_energy_algebraic_5}
        \end{equation}
        again because $\moperator{U}$ and $\moperator{Q}$ are symmetric matrices.

        Finally, substituting (50a) into \eqref{eq:time_derivative_kinetic_energy_algebraic_5} yields
        \begin{equation}
            \ddt{\mathcal{E}_{h}} = -\valgebra{F}^{\transpose}\moperator{U}^{q}\valgebra{F} +  \valgebra{F}^{\transpose}\left(\moperator{E}^{2,1}\right)^{\transpose}\moperator{Q}(\valgebra{K} + g\valgebra{h}) - \valgebra{F}^{\transpose}\left(\moperator{E}^{2,1}\right)^{\transpose}\moperator{Q}(\valgebra{K} + g\valgebra{h}) = 0\,,\label{eq:time_derivative_kinetic_energy_algebraic_6}
        \end{equation}
        because the last two terms directly cancel each other, and the first one is zero since $\moperator{U}^{q}$ is a skew-symmetric matrix.
}
\secondRev{The necessary conditions for energy conservation in space, namely that $\moperator{U}^{q}$ is skew-symmetric and that the 
gradient and divergence operators are anti-adjoints, as given in \eqref{adjoint_divergence} and applied in \eqref{mom_cont_weak_discrete}
are more fundamentally derived from the structure of the Poisson bracket used to construct the rotating shallow water equations in Hamiltonian 
form \cite{Salmon04,ER17}}.
Note that \eqref{eq:time_derivative_kinetic_energy_algebraic_2} requires that the chain rule holds in the discrete form for $\vec u_h$ with 
respect to $t$. The conservation of energy therefore is limited to the truncation error in the time stepping scheme.

\subsection{Conservation of potential enstrophy}

\firstRev{
        Potential enstrophy, $\mathcal{Q}$, is defined as, [4],
        \begin{equation}
            \mathcal{Q} := \int_{\Omega}h q^{2}\,\mathrm{d}\Omega = \innerproduct{hq}{q}\,.
        \end{equation}
        Conservation of potential enstrophy states that
        \begin{equation}\label{pot_enst_time_deriv_expanded}
            \ddt{\mathcal{Q}} = \ddt{}\innerproduct{hq}{q} = \innerproduct{\pddt{h}q}{q} + 2\innerproduct{hq}{\pddt{q}} = 0\,.
        \end{equation}
        When considering the discrete variables, the time variation of potential enstrophy is expressed in an analogous way
        \begin{equation}
            \ddt{\mathcal{Q}_{h}} = \innerproduct{\pddt{h_{h}}q_{h}}{q_{h}} + 2\innerproduct{h_{h}q_{h}}{\pddt{q_{h}}}\,, 
        \end{equation}
        which at the algebraic level becomes
        \begin{equation}
            \ddt{\mathcal{Q}_{h}} = \ddt{\valgebra{h}^{\transpose}}\moperator{W}^{q}\valgebra{q} + 2 \valgebra{h}^{\transpose}\moperator{W}^{q}\ddt{\valgebra{q}}\,, \label{eq:enstrophy_conservation_algebraic_1}
        \end{equation}
        with
        \begin{equation}
            \moperator{W}^{q}_{ij} := \innerproduct{\epsilon^{W}_{i}}{q_{h}\epsilon^{Q}_{j}} \,.
        \end{equation}
        Noting that
        \begin{equation}
            \left(\moperator{W}^{q}\valgebra{h}\right)_{i} = \sum_{k=1}^{d_{Q}}\moperator{W}^{q}_{ik}h_{k} = \sum_{k=1}^{d_{Q}}\innerproduct{\epsilon^{W}_{i}}{q_{h}\epsilon^{Q}_{k}}h_{k} = \sum_{k=1}^{d_{W}}\innerproduct{\epsilon^{W}_{i}}{h_{h}\epsilon^{W}_{k}}q_{k} = \sum_{k=1}^{d_{W}}\moperator{W}^{h}_{ik}q_{k}  = \left(\moperator{W}^{h}\valgebra{q}\right)_{i} \,,
        \end{equation}
        it is possible to rewrite \eqref{eq:enstrophy_conservation_algebraic_1} as
        \begin{equation}
            \ddt{\mathcal{Q}_{h}} = \ddt{\valgebra{h}^{\transpose}}\left(\moperator{W}^{q}\right)^{\transpose}\valgebra{q} + 2 \valgebra{q}^{\transpose}\moperator{W}^{h}\ddt{\valgebra{q}}\,. \label{eq:enstrophy_conservation_algebraic_2}
        \end{equation}
        Furthermore, if we now observe that
        \begin{equation}
            \ddt{}\left(\moperator{W}^{h}\valgebra{q}\right) = \ddt{\moperator{W}^{h}}\valgebra{q} + \moperator{W}^{h}\ddt{\valgebra{q}} = \moperator{W}^{q}\ddt{\valgebra{h}} + \moperator{W}^{h}\ddt{\valgebra{q}} \,,
        \end{equation}
             equation \eqref{eq:enstrophy_conservation_algebraic_2} becomes
             \begin{equation}
            \ddt{\mathcal{Q}_{h}} = \ddt{\valgebra{h}^{\transpose}}\left(\moperator{W}^{q}\right)^{\transpose}\valgebra{q} - 2 \ddt{\valgebra{h}^{\transpose}}\left(\moperator{W}^{q}\right)^{\transpose}\valgebra{q} + 2\valgebra{q}^{\transpose}\ddt{}\left(\moperator{W}^{h}\valgebra{q}\right) = -\ddt{\valgebra{h}^{\transpose}}\left(\moperator{W}^{q}\right)^{\transpose}\valgebra{q} + 2\valgebra{q}^{\transpose}\ddt{}\left(\moperator{W}^{h}\valgebra{q}\right)\,.\label{eq:enstrophy_conservation_algebraic_3}
        \end{equation}
        If we now assume that $\ddt{\valgebra{f}} = 0$, using (50c) we can rewrite \eqref{eq:enstrophy_conservation_algebraic_3} as
        \begin{equation}
            \ddt{\mathcal{Q}_{h}} = -\ddt{\valgebra{h}^{\transpose}}\left(\moperator{W}^{q}\right)^{\transpose}\valgebra{q} - 2\valgebra{q}^{\transpose}\left(\moperator{E}^{1,0}\right)^{\transpose}\moperator{U}\ddt{\valgebra{u}}\,.\label{eq:enstrophy_conservation_algebraic_4}
        \end{equation}
        Substituting (50a) into the last term on the right hand side yields
        \begin{equation}
            \ddt{\mathcal{Q}_{h}} = -\ddt{\valgebra{h}^{\transpose}}\left(\moperator{W}^{q}\right)^{\transpose}\valgebra{q} + 2\valgebra{q}^{\transpose}\left(\moperator{E}^{1,0}\right)^{\transpose}\moperator{U}^{q}\valgebra{F} -  2\valgebra{q}^{\transpose}\left(\moperator{E}^{1,0}\right)^{\transpose}\left(\moperator{E}^{2,1}\right)^{\transpose}\moperator{Q}\left(\valgebra{K} + g\valgebra{h}\right) \,.\label{eq:enstrophy_conservation_algebraic_5}
        \end{equation}
        The last term on the right hand side cancels out because $\moperator{E}^{2,1}\moperator{E}^{1,0} = \moperator{0}$, therefore
        \begin{equation}
            \ddt{\mathcal{Q}_{h}} = -\ddt{\valgebra{h}^{\transpose}}\left(\moperator{W}^{q}\right)^{\transpose}\valgebra{q} + 2\valgebra{q}^{\transpose}\left(\moperator{E}^{1,0}\right)^{\transpose}\moperator{U}^{q}\valgebra{F} \,.\label{eq:enstrophy_conservation_algebraic_6}
        \end{equation}
        Since $\ddt{\mathcal{Q}_{h}}$ is a scalar, we have that
        \begin{equation}
            \ddt{\valgebra{h}^{\transpose}}\left(\moperator{W}^{q}\right)^{\transpose}\valgebra{q} = \valgebra{q}^{\transpose}\moperator{W}^{q}\ddt{\valgebra{h}}\,,
        \end{equation}
        and substituting into \eqref{eq:enstrophy_conservation_algebraic_6}, results in
        \begin{equation}
            \ddt{\mathcal{Q}_{h}} = -\valgebra{q}^{\transpose}\moperator{W}^{q}\ddt{\valgebra{h}}+ 2\valgebra{q}^{\transpose}\left(\moperator{E}^{1,0}\right)^{\transpose}\moperator{U}^{q}\valgebra{F} \,.\label{eq:enstrophy_conservation_algebraic_7}
        \end{equation}
        Moreover, using (50b) we can write
        \begin{equation}
            \ddt{\mathcal{Q}_{h}} = -\valgebra{q}^{\transpose}\moperator{W}^{q}\moperator{E}^{2,1}\valgebra{F}+ 2\valgebra{q}^{\transpose}\left(\moperator{E}^{1,0}\right)^{\transpose}\moperator{U}^{q}\valgebra{F} \,.\label{eq:enstrophy_conservation_algebraic_8}
        \end{equation}
        Now, noting that
        \begin{equation}
            \left[\left(\moperator{U}^{q}\right)^{\transpose}\moperator{E}^{1,0}\valgebra{q}\right]_{j} = \sum_{i,k=1}^{d_{W},d_{U}}\innerproduct{\vfield{\epsilon}^{\,U}_{k}}{q_{h}\times\vfield{\epsilon}^{\,U}_{j}}\moperator{E}^{1,0}_{k,i}q_{i} = - \sum_{i,k=1}^{d_{W},d_{U}}\innerproduct{\vfield{\epsilon}^{\,U}_{j}}{q_{h}\times\vfield{\epsilon}^{\,U}_{k}}\moperator{E}^{1,0}_{k,i}q_{i} = - \innerproduct{\vfield{\epsilon}^{\,U}_{j}}{q_{h}\times\left(\sum_{i,k=1}^{d_{W},d_{U}}\moperator{E}^{1,0}_{k,i}q_{i}\vfield{\epsilon}^{\,U}_{k}\right)} \,. \label{eq:enstrophy_conservation_algebraic_9}
        \end{equation}
        Recalling (42), we can rewrite the last term of \eqref{eq:enstrophy_conservation_algebraic_9} as
        \begin{equation}
            \sum_{i,k=1}^{d_{W},d_{U}}\moperator{E}^{1,0}_{k,i}q_{i}\vfield{\epsilon}^{\,U}_{k} = \nabla^{\perp}q_{h}\,,
        \end{equation}
         and consequently \eqref{eq:enstrophy_conservation_algebraic_9} becomes
        \begin{equation}
            \left[\left(\moperator{U}^{q}\right)^{\transpose}\moperator{E}^{1,0}\valgebra{q}\right]_{j}  =  - \innerproduct{\vfield{\epsilon}^{\,U}_{j}}{q_{h}\times\nabla^{\perp}q_{h}}  \,. \label{eq:enstrophy_conservation_algebraic_10a}
        \end{equation}
        The right hand side of this equation can be transformed in the following way
        \begin{equation}
            - \innerproduct{\vfield{\epsilon}^{\,U}_{j}}{q_{h}\times\nabla^{\perp}q_{h}} \stackrel{*}{=} \frac{1}{2}\innerproduct{\vfield{\epsilon}^{\,U}_{j}}{\nabla(q_{h}\cdot{q_{h}})} \stackrel{*}{=} \frac{1}{2}\innerproduct{\nabla\cdot\vfield{\epsilon}^{\,U}_{j}}{q_{h}\cdot{q_{h}}}\,. \label{eq:enstrophy_conservation_algebraic_10}
        \end{equation}
        Therefore we have
        \begin{equation}
            \left[\left(\moperator{U}^{q}\right)^{\transpose}\moperator{E}^{1,0}\valgebra{q}\right]_{j} = \frac{1}{2}\innerproduct{\nabla\cdot\vfield{\epsilon}^{\,U}_{j}}{q_{h}\cdot{q_{h}}}\,. \label{eq:enstrophy_conservation_algebraic_11}
        \end{equation}
        It is important to note that the two middle identities, with $\stackrel{*}{=}$, in \eqref{eq:enstrophy_conservation_algebraic_10} 
are only valid if exact integration is used in the inner products. If approximate integration is used, these identities are only guaranteed 
to be approximately valid. Now, the right hand term in \eqref{eq:enstrophy_conservation_algebraic_11} can be written as
        \begin{equation}
            \innerproduct{\nabla\cdot\vfield{\epsilon}^{\,U}_{j}}{q_{h}\cdot{q_{h}}} = \sum_{k=1}^{d_{Q}}\moperator{E}^{2,1}_{k,j}\innerproduct{\epsilon^{Q}_{k}}{q_{h}\cdot{q_{h}}} =  \sum_{k=1}^{d_{Q}}\moperator{E}^{2,1}_{k,j}\innerproduct{\epsilon^{Q}_{k}q_{h}}{\epsilon^{W}_{i}}q_{i} = \left[\left(\moperator{E}^{2,1}\right)^{\transpose}\left(\moperator{W}^{q}\right)^{\transpose}\valgebra{q}\right]_{j}\,.
        \end{equation}
        Therefore, \eqref{eq:enstrophy_conservation_algebraic_8} becomes
        \begin{equation}
            \ddt{\mathcal{Q}_{h}} = -\valgebra{q}^{\transpose}\moperator{W}^{q}\moperator{E}^{2,1}\valgebra{F}+ \valgebra{q}^{\transpose}\moperator{W}^{q}\moperator{E}^{2,1}\valgebra{F} = 0 \,,
        \end{equation}
        thus proving conservation of potential enstrophy at the discrete level.
}

As for energy conservation, potential enstrophy conservation requires that the chain rule holds for differentiation in time, as applied in
\eqref{pot_enst_time_deriv_expanded}. However unlike the energy conservation argument, potential enstrophy conservation
also requires that this holds for spatial differentiation, as given in \eqref{eq:enstrophy_conservation_algebraic_10}. Potential enstrophy
conservation is therefore also subject to exact spatial integration in the weak form.

\subsection{Geostrophic balance}

The existence of stationary geostrophic modes derives from the linearization of \eqref{eq:shallow_water_equations_all} as \cite{CS12}

\begin{subnumcases}{\label{eq:shallow_water_equations_linear}}
    \frac{\partial\vec u}{\partial t} + f\times\vec u + g\nabla h = 0\,, & in $\Omega\times (0,t_{F}]\,,$  \label{mom_cont_linear}\\
    \frac{\partial h}{\partial t} + H\nabla\cdot\vec u = 0\,, & in $\Omega\times (0,t_{F}]\,,$ \label{mas_cont_linear}
\end{subnumcases}
where $H$ is the mean depth of the fluid layer.
If we assume a stationary solution then we have the balanced system

\begin{equation}\label{geostrophic_balance_1}
    f\times\vec u + g\nabla h = 0\qquad\nabla\cdot\vec u = 0
\end{equation}
In order to show that this balance is satisfied for our formulation we begin by introducing the stream function
as $\vec u = \nabla^{\perp}\psi$. Note that \secondRev{this } identity is satisfied in the strong form for our discrete formulation
via \eqref{eq:hilbert_subcomplex_basis_W}. Substituting this into \eqref{geostrophic_balance_1} gives

\begin{equation}\label{geostrophic_balance_2}
    -f\nabla\psi + g\nabla h = 0\qquad\nabla\cdot\vec u = 0
\end{equation}
To show that this relation is satisfied at the algebraic level we begin by taking the linearized form of \eqref{mom_cont_weak_discrete_matrix},
\eqref{mas_cont_weak_discrete_matrix} as

\firstRev{
\begin{equation}\label{mom_cont_weak_discrete_linear}
\boldsymbol{\mathsf{U}}\frac{\mathrm{d}\boldsymbol{u}}{\mathrm{d}t} + \boldsymbol{\mathsf{U}}^{f}\boldsymbol{u} - 
g\left(\boldsymbol{\mathsf{E}}^{2,1}\right)^{\top}\boldsymbol{\mathsf{Q}}\boldsymbol{h} = 0
\end{equation}
\begin{equation}\label{mas_cont_weak_discrete_linear}
\frac{\mathrm{d}\boldsymbol{h}}{\mathrm{d}t} + \boldsymbol{\mathsf{E}}^{2,1}\boldsymbol{u} = 0
\end{equation}
where $\mathsf{U}^{f}_{ij} := \ip{f_{h}\times\vec{\epsilon}_{j}^{\,U}}{\vec{\epsilon}_{i}^{\,U}}$.
Note that we have eliminated the $\boldsymbol{\mathsf{Q}}$ matrix in the linearized continuity equation 
\eqref{mas_cont_weak_discrete_linear} as done in \eqref{eq:time_variation_of_mass_and_flux}. 
}
Applying the stream function identity 
to \eqref{mom_cont_weak_discrete_linear} and then the adjoint relation \eqref{adjoint_divergence} gives
\firstRev{
\begin{equation}\label{mom_cont_weak_discrete_linear_2}
\boldsymbol{\mathsf{U}}\frac{\mathrm{d}\boldsymbol{u}}{\mathrm{d}t} +
f\left(\boldsymbol{\mathsf{E}}^{2,1}\right)^{\top}\boldsymbol{\mathsf{Q}}^w\boldsymbol{\psi} -
g\left(\boldsymbol{\mathsf{E}}^{2,1}\right)^{\top}\boldsymbol{\mathsf{Q}}\boldsymbol{h} = 0
\end{equation}
where $\mathsf{Q}^{w}_{ij} := \ip{{\epsilon}_{j}^{\,W}}{{\epsilon}_{i}^{\,Q}}$.
}
The last two terms of \eqref{mom_cont_weak_discrete_linear_2} constitute the weak form of \eqref{geostrophic_balance_2} and so 
\firstRev{balance up to truncation error of the interpolating polynomials, $\epsilon^W_i$, $\epsilon^Q_i$ is achieved }
such that \firstRev{$\boldsymbol{u }$ is approximately } stationary. Since \eqref{mas_cont_weak_discrete_linear} 
holds point wise, the absence of divergence leads to a stationary fluid depth \firstRev{$\boldsymbol{h}$}, 
thus demonstrating geostrophic balance.

In the following section we will illustrate for a simple divergence free linear test case that the errors
in geostrophic balance remain bounded and convergent with temporal and spatial resolution.

\section{Results}

\subsection{Convergence}

In order to first validate our mixed mimetic spectral element formulation we inspect the $L_2$ norm error convergence
of the various operators. For this we compare the three diagnostic equations
(\ref{eq:definition_potential_vorticity_weak_discrete_matrix}-\ref{eq:definition_kinetic_energy_weak_discrete_matrix})
to the analytic solutions for a specified velocity and depth field, as derived from a stream function solution of

\begin{equation}\label{strm_func}
\psi = 0.1\cos(x-\pi)\cos(y-\pi)\qquad \Omega = (0,2\pi]\times (0,2\pi],
\end{equation}
where $\vec u = \nabla^{\perp}\psi$, $h = (f/g)\psi + H$ via geostrophic balance, $f = g = 8.0$ and
$H = 0.2$ is the constant of integration in the geostrophic balance relation. In each case exact spatial
integration is applied.
By demonstrating both the algebraic convergence of the errors for constant
polynomial degree with decreasing mesh size and the spectral convergence with increasing polynomial degree and
constant mesh size of the solution for each of the diagnostic equations we aim to
validate the basis functions for each of the $W_h$, $U_h$ and $Q_h$ function spaces.

\begin{figure}[!hbtp]
\includegraphics[width=0.48\textwidth,height=0.36\textwidth]{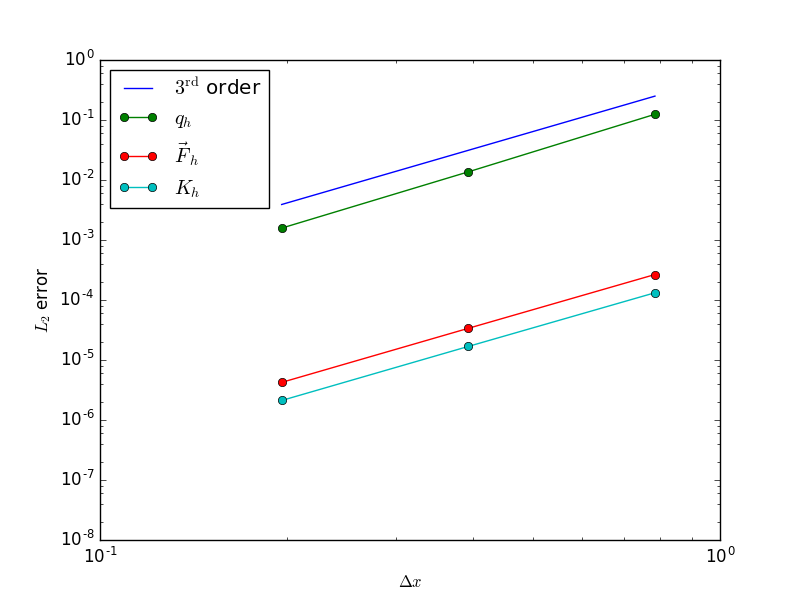}
\includegraphics[width=0.48\textwidth,height=0.36\textwidth]{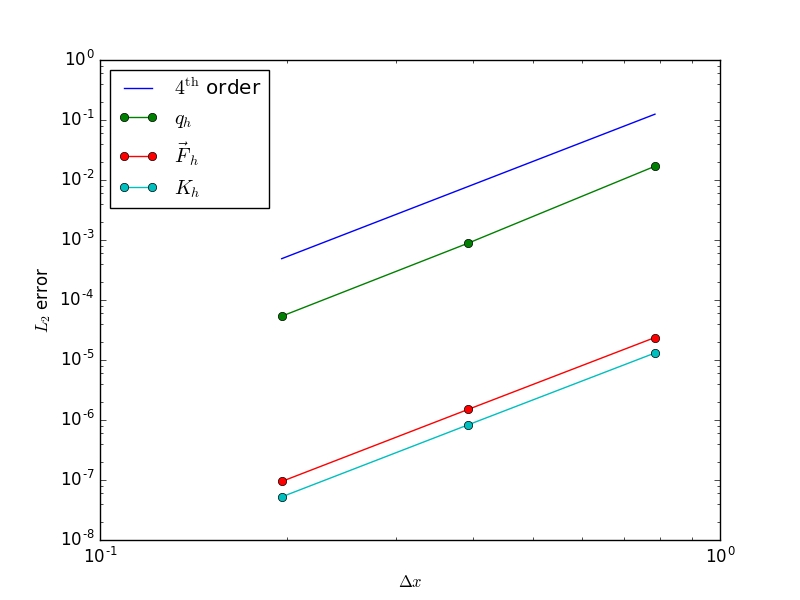}\\
\caption{Algebraic convergence for the $3^{\mathrm{rd}}$ order (left) and $4^{\mathrm{th}}$ order
(right) basis functions for $q_h$, $\vec F_h$ and $K_h$.
Blue slopes show the theoretical convergence rates for $3^{\mathrm{rd}}$ and $4^{\mathrm{th}}$
order respectively.}
\label{fig::diagnostic_convergence}
\end{figure}
As shown in \figref{fig::diagnostic_convergence}, $\vec F_h$ and $K_h$ satisfy their expected 
rates of convergence for both $3^\mathrm{rd}$ and $4^\mathrm{th}$ order accurate basis representations.
However $q_h \in W_h$ should in principle converge at one order higher, since its polynomial expansion
is one degree higher. This higher rate of convergence is not observed from \figref{fig::diagnostic_convergence}.
The fact that $q_h$ converges at the same rate as $\vec F_k$ and $K_h$ may be possibly be attributed
to the fact that $h_h$ as used in the left hand side of diagnostic equation \eqref{eq:definition_potential_vorticity_weak_discrete_matrix}
is discontinuous across element boundaries and so may be projecting spurious gradients onto $q_h$.
Convergence studies of $q_h$ using the $||q_h||^2_{H(\mathrm{rot})}$ norm with an analytic (continuous)
depth field (not shown), which exhibit the correct rate of convergence, would seem to suggest this.

\begin{figure}[!hbtp]
\includegraphics[width=0.48\textwidth,height=0.36\textwidth]{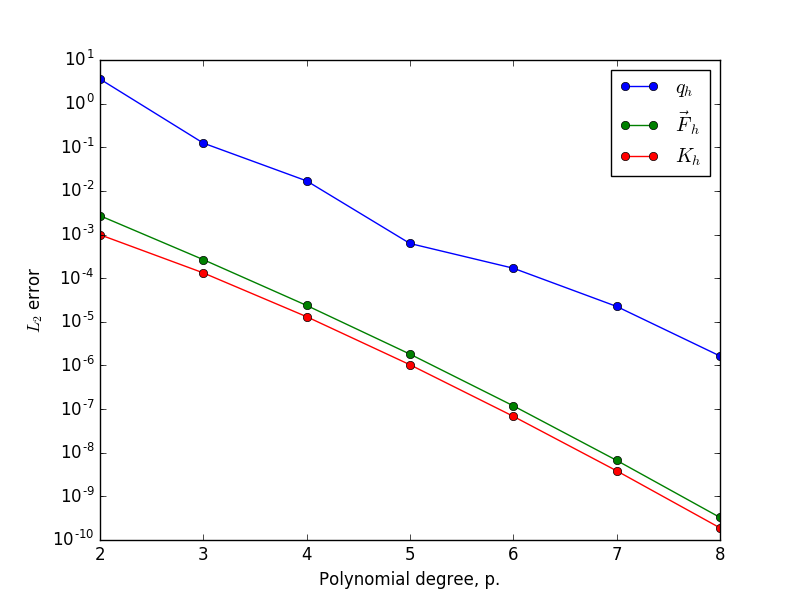}
\includegraphics[width=0.48\textwidth,height=0.36\textwidth]{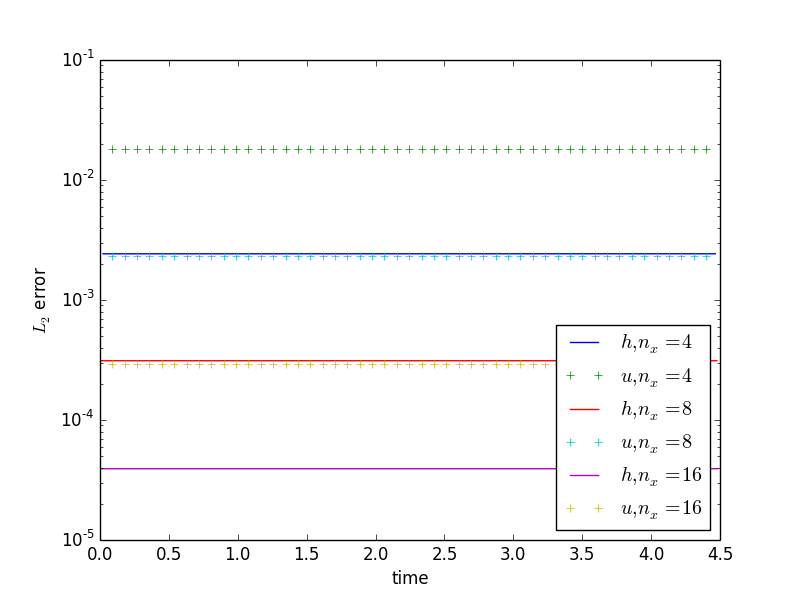}\\
\caption{Spectral convergence of the diagnostic equations with polynomial degree (left).
and convergence of errors for the linear geostrophic balance test (right).}
\label{fig::geostrophic_balance}
\end{figure}

\figref{fig::geostrophic_balance} shows the spectral convergence of errors for constant
mesh size and increasing polynomial degree for each of the 
$W_h$, $U_h$ and $Q_h$ tensor product element diagnostic equations.

\subsection{Linearized geostrophic balance}

As a second test the linearized system (\ref{mom_cont_linear},\ref{mas_cont_linear}) is solved in
the discrete form for a divergence free velocity field and a depth field derived from geostrophic
balance via the stream function solution \eqref{strm_func}. The discrete
problem for this system is equivalent to \eqref{eq:shallow_water_continuous_weak_form_discrete_matrix} with
$\boldsymbol{q} = \boldsymbol{f}$, $\boldsymbol{K} = 0$ and $\boldsymbol{F} = H\boldsymbol{u}$.
\firstRev{For this and all subsequent tests we use an explicit second order Runge-Kutta scheme, given 
for the continuous form \eqref{eq:shallow_water_equations_all} as 
\[y' = y^n - 0.5\Delta tG(y^n),\qquad y^{n+1} = y^n - \Delta tG(y')\]
for $y = [\vec u, h]^{\top}$, $G = [q\times\vec F + \nabla(K + gh),\nabla\cdot\vec F]^{\top}$.}
\figref{fig::geostrophic_balance} shows the convergence of errors for the linearized geostrophic balance
problem with $n_x = 4, 8$ and $16$ $3^{\mathrm{rd}}$ order elements in each dimension (and a time step of
$\Delta t = 0.02/n_x$). As can be seen the errors for a given resolution remain constant,
showing that once the weak form geostrophic balance relation has been approximated for the first time
step there is no subsequent accumulation of errors.

\subsection{Conservation}

\begin{figure}[!hbtp]
\includegraphics[width=0.48\textwidth,height=0.36\textwidth]{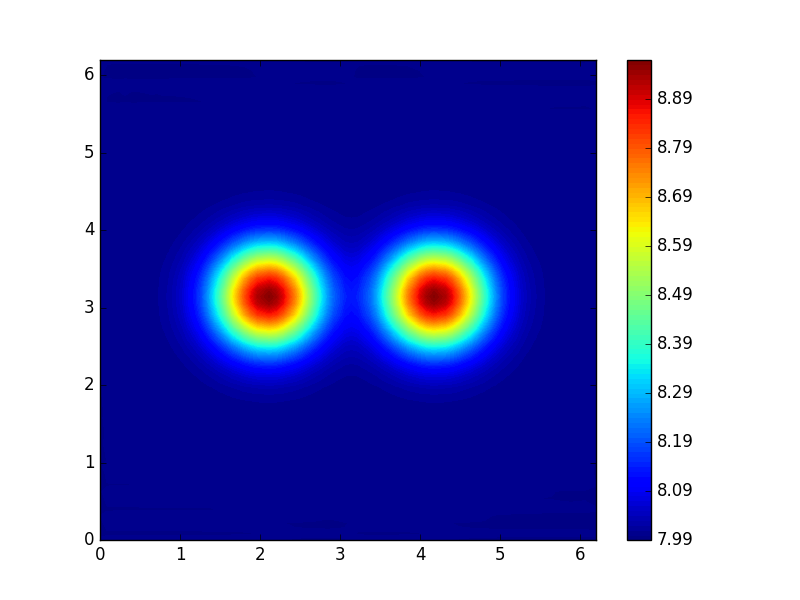}
\includegraphics[width=0.48\textwidth,height=0.36\textwidth]{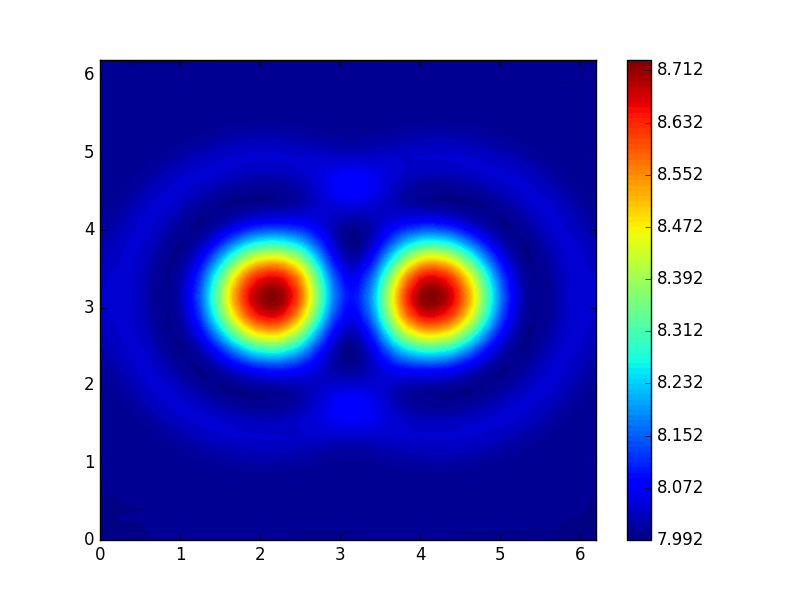}\\
\includegraphics[width=0.48\textwidth,height=0.36\textwidth]{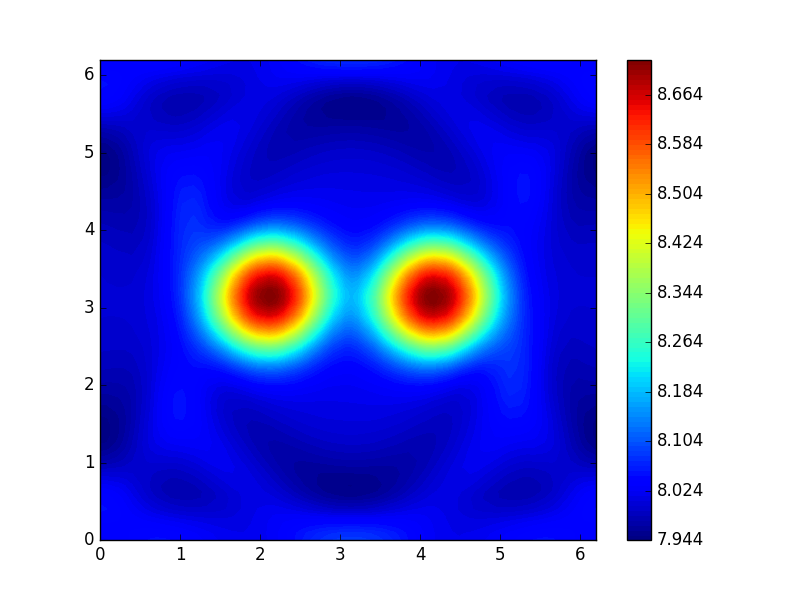}
\includegraphics[width=0.48\textwidth,height=0.36\textwidth]{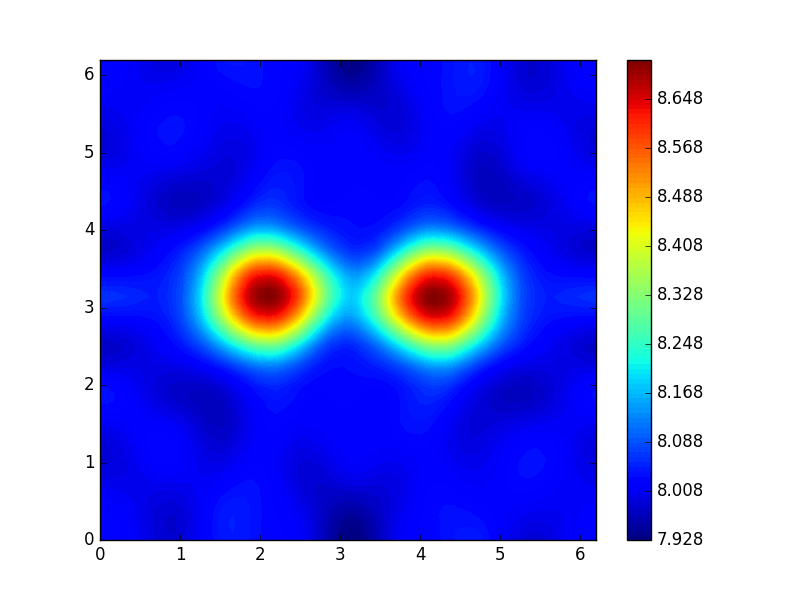}\\
\caption{Fluid depth $h$ at times $t = $ $0.0$, $0.5$, $1.0$ and $2.0$
dimensionless units. Approximate geostrophic balance is preserved while fast gravity waves
radiate outwards from the disturbances.}
\label{fig::conservation_exact}
\end{figure}

Having verified the construction of our mixed mimetic spectral element operators and solver via
various convergence and geostrophic balance tests, we proceed to investigate the conservation properties
of the full system \eqref{eq:shallow_water_continuous_weak_form_discrete_matrix}, as derived in Section 4.
The model is tested for a pair of vortices starting in approximate geostrophic balance, as given
from the stream function solution 

\begin{equation}
\psi = e^{-2.5((x-\pi)^2 + (y-2\pi/3)^2)} + e^{-2.5((x-\pi)^2 + (y-4\pi/3)^2)}\qquad \Omega = (0,2\pi]\times (0,2\pi],
\end{equation}
where the velocity is given as $\vec u = \nabla^{\perp}\psi$, and the depth is derived from geostrophic balance as
$f\times\vec u + g\nabla h = 0$, with $f = g = H = 8.0$, using $20\times 20$ $3^{\mathrm{rd}}$ order elements
and a maximum time step of $\Delta t = 0.0052$. The solution behaves as expected, with
fast gravity waves radiating out from the initial disturbance, which is preserved for long times
due to the approximate geostrophic balance of the initial condition, as shown in \figref{fig::conservation_exact}.
The ratio of the deformation radius, $L_d = \sqrt{gH}/f$ to the nodal grid spacing is approximately 9.55.

\begin{figure}[!hbtp]
\includegraphics[width=0.48\textwidth,height=0.42\textwidth]{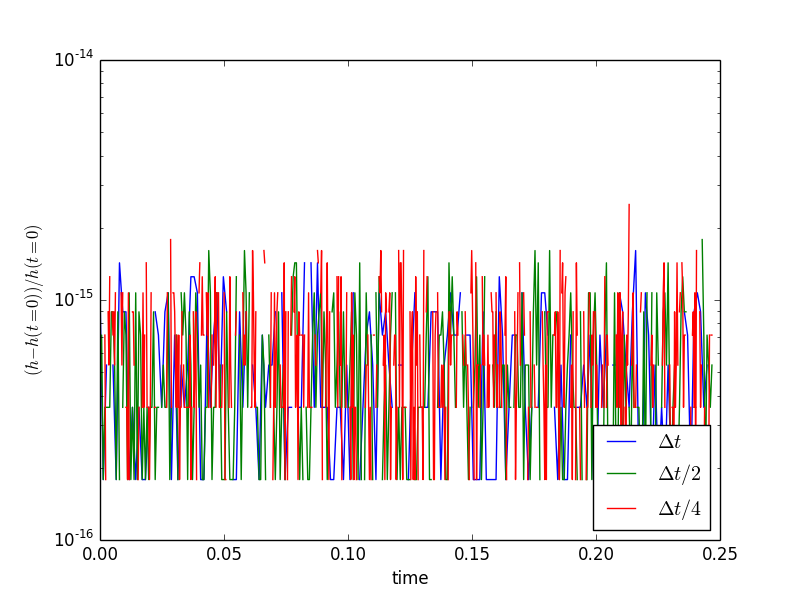}
\includegraphics[width=0.48\textwidth,height=0.42\textwidth]{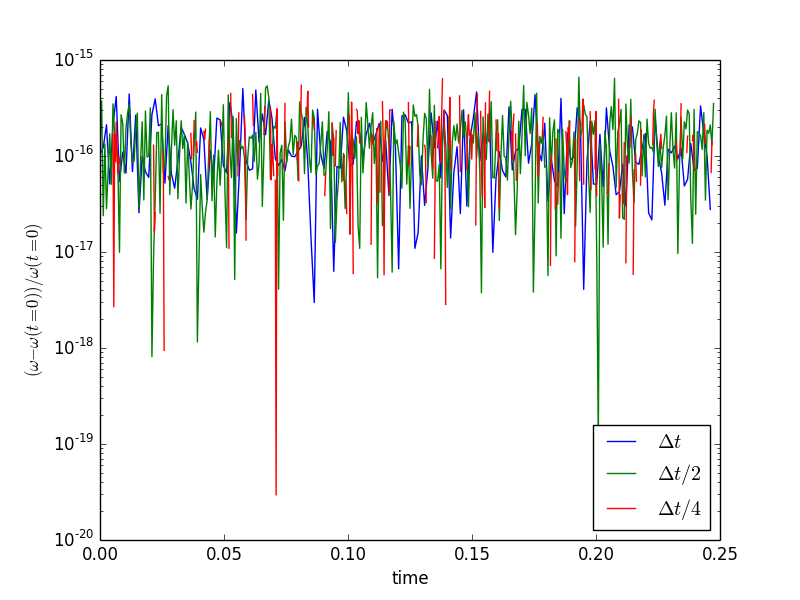}\\
\caption{Exact conservation for the volume, $h$ (left) and vorticity, $\omega$ (right) with time.}
\label{fig::conservation_hw}
\end{figure}

\begin{figure}[!hbtp]
\includegraphics[width=0.48\textwidth,height=0.42\textwidth]{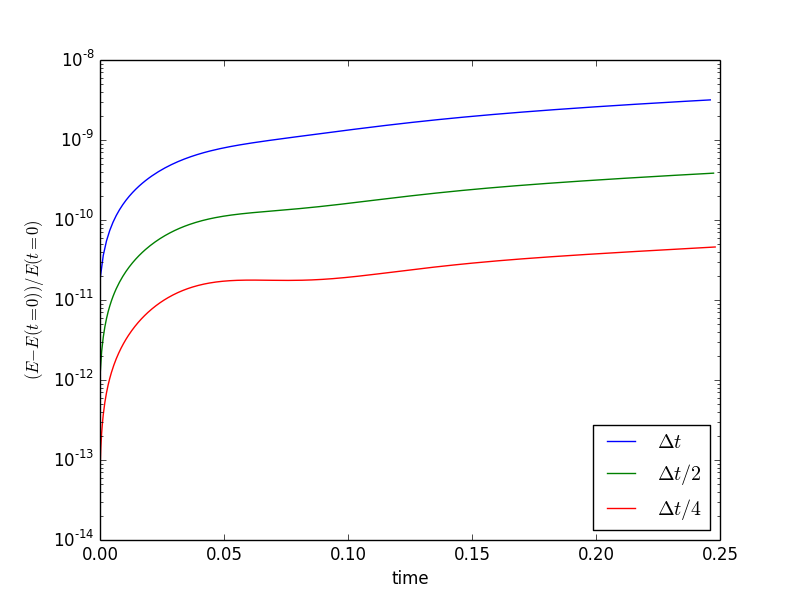}
\includegraphics[width=0.48\textwidth,height=0.42\textwidth]{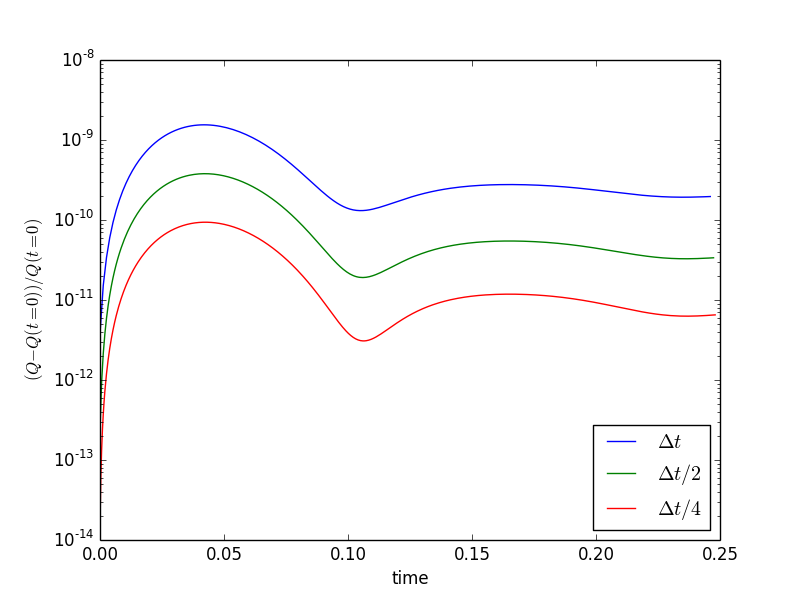}\\
\caption{Convergence of conservation errors for the energy, $\mathcal{E}$ (left) and potential enstrophy, $\mathcal{Q}$ (right)
with time step $\Delta t$ (exact spatial integration).}
\label{fig::conservation_EQ_exact}
\end{figure}

As discussed in Section 4, mass and vorticity conservation hold independent of time step, as shown in
\figref{fig::conservation_hw}. This is due to the point wise satisfaction of the divergence theorem in
the case of mass, and the elimination of the gradient operator by the curl in the weak form
in the case of vorticity. Total energy and potential enstrophy are conserved to truncation
error in time, as shown in \figref{fig::conservation_EQ_exact}, with the second order Runge-Kutta time stepping scheme
applied with varying time step size.

\begin{figure}[!hbtp]
\includegraphics[width=0.48\textwidth,height=0.42\textwidth]{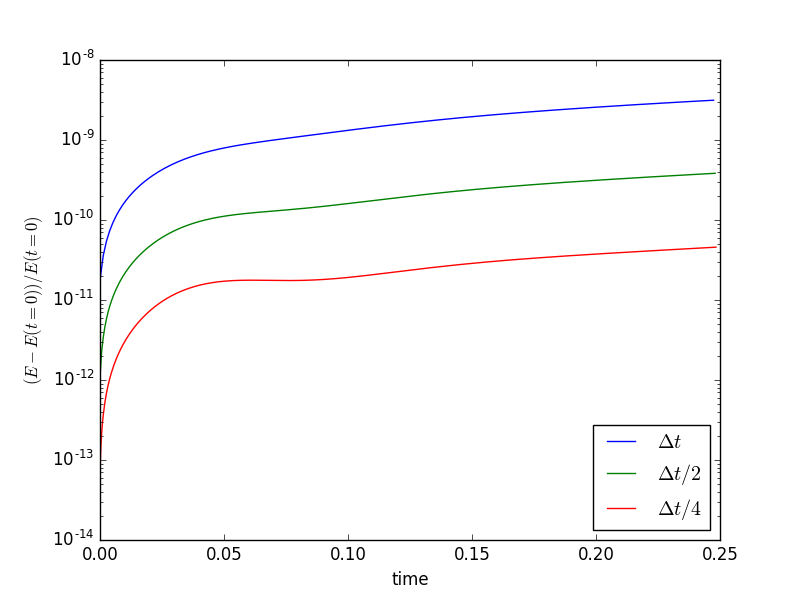}
\includegraphics[width=0.48\textwidth,height=0.42\textwidth]{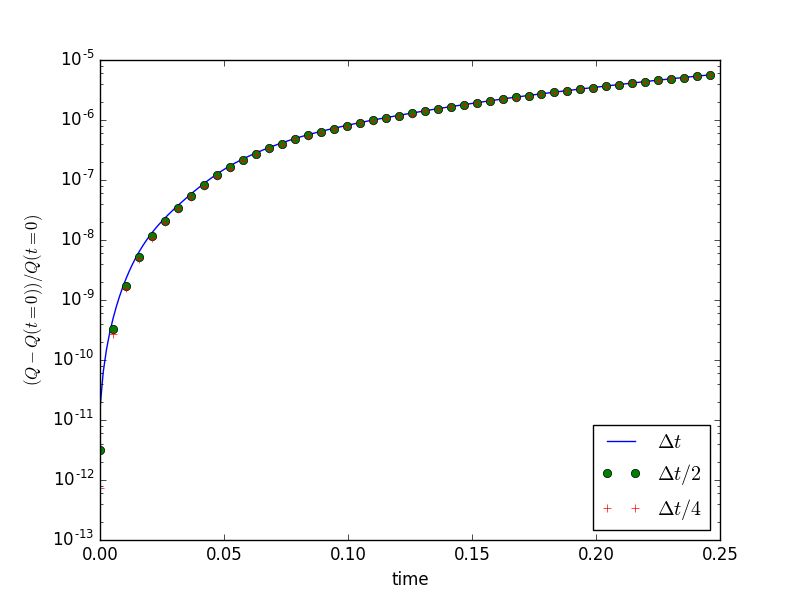}\\
\caption{Convergence of conservation errors for the energy, $\mathcal{E}$ (left) and potential enstrophy, $\mathcal{Q}$ (right)
with time step $\Delta t$ (inexact spatial integration).}
\label{fig::conservation_EQ_inexact}
\end{figure}

\secondRev{Gauss-Lobatto-Legendre (GLL) quadrature is known to be exact for polynomials of degree $p = 2n - 3$, where 
$n$ is the number of quadrature points \cite{KS05}. Therefore in order to exactly integrate all nonlinear matrices in 
\eqref{eq:shallow_water_continuous_weak_form_discrete_matrix} we use $n = (3p + 3)/2$ quadrature points (where $3p$ 
is the maximum degree of the test function, trial function and nonlinear function basis expansion product). 
A second set of tests was also run using inexact quadrature with $n = p + 1$ in order to derive diagonal mass matrices 
for test functions in $W_h$, since this significantly increases the computational efficiency of the scheme, and is 
customary for spectral element models \cite{TF10}.}
While energy conservation holds for inexact spatial integration, 
potential enstrophy conservation fails for inexact integration, as shown in \figref{fig::conservation_EQ_inexact}.


\subsection{Stabilization of shear flow via removal of anticipated potential vorticity}

Our final numerical experiment involves an investigation of the anticipated potential vorticity method
\cite{SB85} for a shear flow initial condition over a stationary orographic feature.
The orography is implemented as a $b_h \in Q_h$, such that the momentum equation \eqref{mom_cont_weak_discrete_matrix} becomes

\begin{equation}
\boldsymbol{\mathsf{U}}\frac{\mathrm{d}\boldsymbol{u}}{\mathrm{d}t} + \boldsymbol{\mathsf{U}}^{q}\boldsymbol{F} -
\left(\boldsymbol{\mathsf{E}}^{2,1}\right)^{\top}\boldsymbol{\mathsf{Q}}\,(\boldsymbol{K}+g\boldsymbol{h} + g\boldsymbol{b}) = 0
\end{equation}
and the energy is correspondingly defined as
$\mathcal{E} = \ip{K_h + 0.5gh_h + gb_h}{h_h}$.
The initial conditions are given as

\begin{equation}
h = H + 0.1\tanh\Big(\frac{1 - y^2}{2}\Big)\qquad
\vec u = \Big(-\frac{\partial h}{\partial y}, 0\Big)
\end{equation}
and the orography as

\begin{equation}
b = \left \{ \begin{array}{ll}
0.0125(\cos(4\pi x/L) + 1)(\cos(4\pi y/L) + 1) & \mbox{if } |x| \le L/4, |y| \le L/4\\[1ex]
0 & \mathrm{otherwise}
\end{array} \right .\;
\end{equation}
with $\Omega = (-L/2,+L/2]\times (-L/2,L/2]$, $L = 10$ and $f = g = H = 1.0$. In order to stabilize the potential enstrophy cascade,
the anticipated potential vorticity $\hat q_h$ \eqref{apv} substitutes for the potential vorticity $q_h$ in the
rotational term of the momentum equation. This results in a loss of potential enstrophy such that the cascade to
sub grid scales is arrested \cite{SB85}. The test is run on $24\times 24$ $3^{\mathrm{rd}}$ order elements, with
a ratio of $L_d$ to the average nodal grid spacing of 7.2.

\begin{figure}[!hbtp]
\includegraphics[width=0.48\textwidth,height=0.36\textwidth]{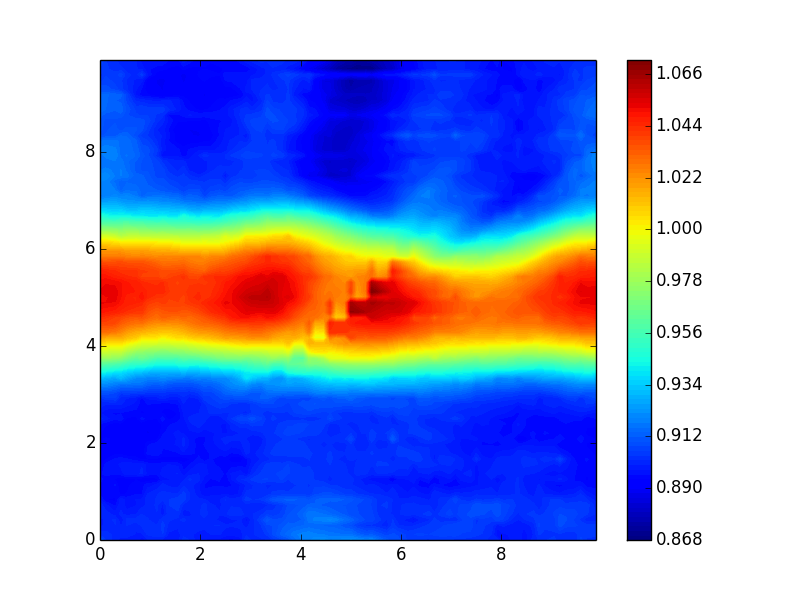}
\includegraphics[width=0.48\textwidth,height=0.36\textwidth]{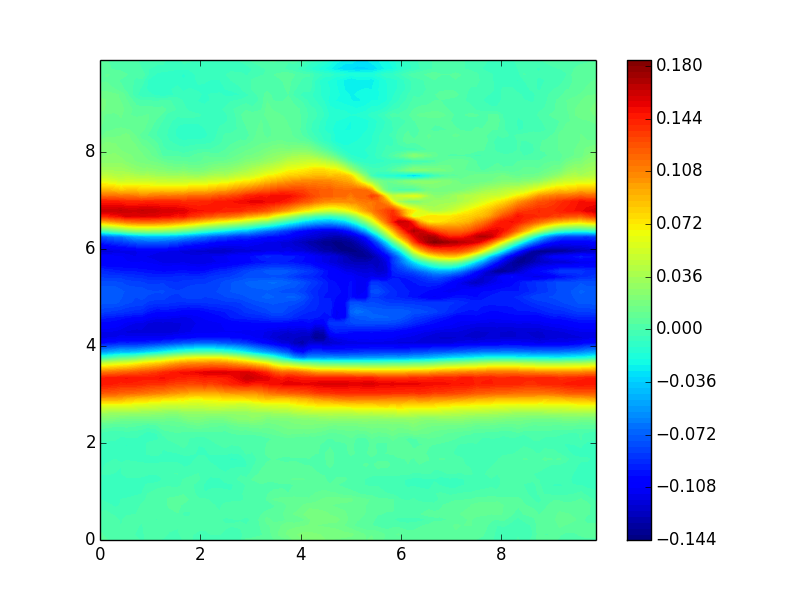}\\
\caption{Fluid depth $h$ (left) and vorticity $\omega$ (right) fields for shear flow over
orography at $t=44$, $\Delta\tau = 0.02$.}
\label{fig::sea_mount_1}
\end{figure}
\begin{figure}[!hbtp]
\includegraphics[width=0.48\textwidth,height=0.36\textwidth]{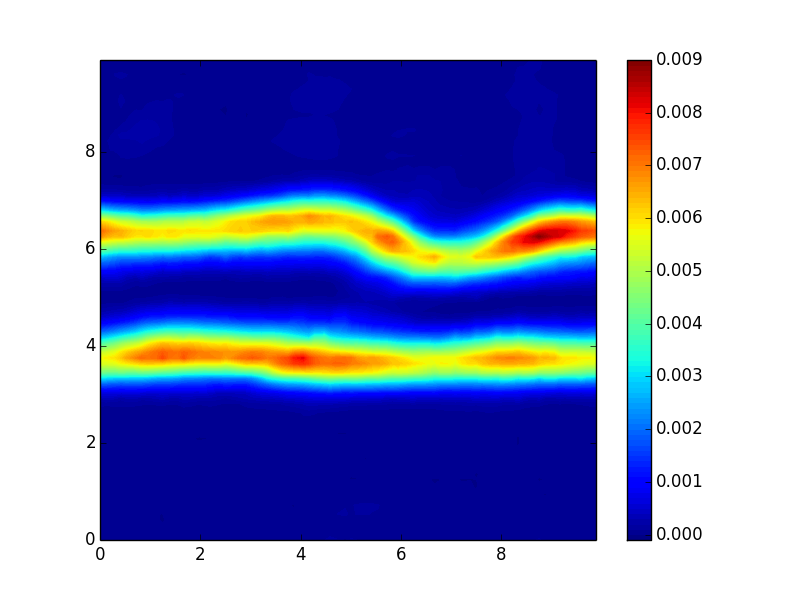}
\includegraphics[width=0.48\textwidth,height=0.36\textwidth]{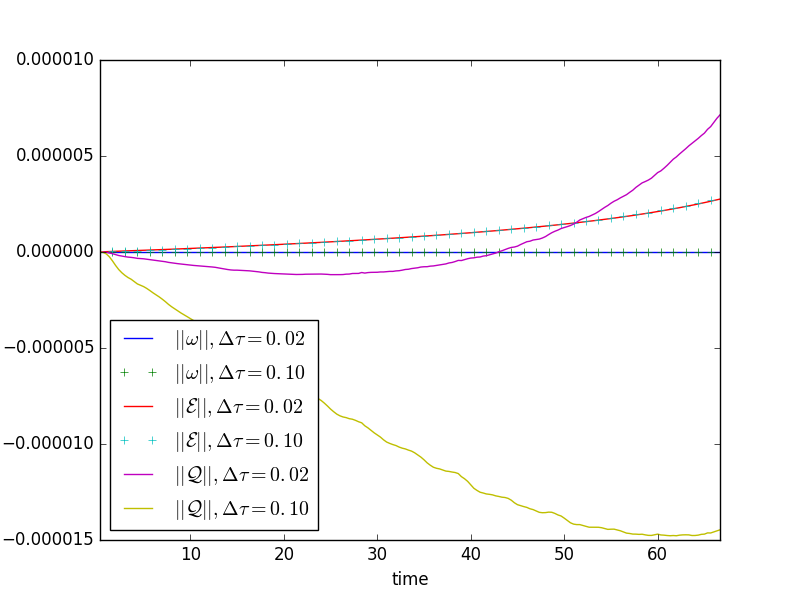}\\
\caption{Kinetic energy per unit volume $K$ for shear flow over orography at $t=44$,
$\Delta\tau = 0.02$ (left), and growth of normalized vorticity $\omega$, energy $\mathcal{E}$ and
potential enstrophy $\mathcal{Q}$ with time for varying values of the anticipated potential vorticity
coefficient $\Delta\tau$, where $||A|| = (A(t) - A(t=0))/A(t=0)$.}
\label{fig::sea_mount_2}
\end{figure}

\figref{fig::sea_mount_1} and \figref{fig::sea_mount_2} show the depth, vorticity and kinetic energy fields
after 44 units of dimensionless time have elapsed. While fast gravity waves radiate out from the topography
as the solution adjusts, the vorticity field evolves on a slow time scale. The outline of the discontinuous
$Q_h$ elements of the orography field is visible in the fluid depth $h_h$. \figref{fig::sea_mount_2} also
shows the growth in vorticity $\omega$, energy $\mathcal{E}$, and potential enstrophy $\mathcal{Q}$ with time
for varying values of the anticipated potential vorticity time scale $\Delta\tau$. As can be seen, vorticity
conservation is preserved to machine precision independent of $\Delta\tau$, and the rate of energy
growth is also similar for both values of $\Delta\tau$. The rates of potential enstrophy growth differ however,
with $\Delta\tau = 0.02$ resulting in the faster growth of potential enstrophy (following an initial loss)
than $\Delta\tau = 0.1$. Neither values of $\Delta\tau$ is able to fully suppress the growth of $\mathcal{Q}$,
which over long time integrations will cascade to sub grid scales and contaminate the solution. This failure
to fully stabilize the solution via the removal of anticipated potential vorticity is perhaps a result of the 
projection of the discontinuous $h_h$ field onto $q_h$ within the potential vorticity diagnostic equation.
Some method of averaging of $h_h$ at the element boundary quadrature points within \eqref{eq:definition_potential_vorticity_weak_discrete_matrix}
may help to ameliorate this.

\section{Conclusion}

In this paper we have built upon the work of previous authors in the development of compatible
finite element methods for geophysical fluid dynamics \cite{CS12,MC14,NSC16} in order to present
a shallow water solver which exactly conserves \secondRev{first order and higher order } moments using the mixed
mimetic spectral element method \cite{Gerritsma11,KG13,Hiemstra14}. The conservation of second
order moments (energy and potential enstrophy) is subject to the truncation error in the time
integration scheme, and the conservation of potential enstrophy also requires exact spatial
integration, as shown by the conservation arguments and demonstrated for the test cases given above.

We note the performance constraints of the different diagnostic and prognostic equations as
follows:
\begin{itemize}
\item Fluid depth, $h_h$: The continuity equation \eqref{mas_cont_weak_discrete_matrix}
is satisfied point wise in the strong form, and may be evaluated purely from the topology, with the
divergence theorem satisfied exactly such that the change in fluid depth is simply the sum of the
momentum fluxes across the adjacent $U_h$ basis functions of $\vec F_h$, and so is extremely fast to compute.

\item Kinetic energy per unit volume, $K_h$: Since $K_h$ also exists on the
discontinuous spaces of $Q_h$, the weak form diagnostic equation \eqref{eq:definition_kinetic_energy_weak_discrete_matrix} may be solved
as a discontinuous Galerkin problem without the need for a global matrix solve.

\item Potential vorticity, $q_h$: If we are prepared to sacrifice potential enstrophy conservation
for an inexact quadrature rule, then the left hand side of \eqref{eq:definition_potential_vorticity_weak_discrete_matrix} is diagonal due
to the orthogonality of the $W_h$ basis functions. This again avoids the need for a
global matrix solve. \secondRev{The use of inexact GLL quadrature also reduces the computational
cost of assembling all other matrices and vectors. This saving may prove significant and override the
desirability of potential enstrophy conservation for production codes.}

\item Velocity, $\vec u_h$ and momentum, $\vec F_h$: Equations \eqref{mom_cont_weak_discrete_matrix} and
\eqref{eq:definition_h_flux_weak_discrete_matrix} require
a global matrix solve, since the function space $U_h$ is continuous across element boundaries.
The solution of these equations therefore represents the major computational bottleneck of the scheme.
\end{itemize}

While we have derived the potential vorticity from a diagnostic equation \eqref{eq:definition_potential_vorticity_weak_discrete_matrix}
in our current formulation, this could alternatively be derived from a prognostic equation by taking
the curl of the momentum equation (ie: by substituting \eqref{mom_cont_weak_discrete_matrix} into \eqref{eq:time_derivative_definition_vorticity}).
Doing so may have the added advantage of allowing for conservation of energy and potential enstrophy
independent of time step via a time staggering of the variables, as has been previously shown for the
2D Navier-Stokes equations \cite{PG17}.

In order to compare the properties of the mixed mimetic spectral element method to the standard
A-grid spectral element method \cite{TF10}, the gravity wave dispersion relation for the two method
will also be compared, in order to determine if the mixed mimetic method improves upon the spurious
discontinuities present in the standard spectral method dispersion relation \cite{MST12}. The power
spectra for nonlinear problems will also be compared to determine if the
mixed method can be run with lower amounts of diffusion due to the absence of collocated velocity
and pressure degrees of freedom.

\section*{Acknowledgements}

David Lee would like to thank Drs. Chris Eldred and Ren\'{e} Hiemstra for their helpful discussions
and insights.
This research was supported as part of the Launching an Exascale ACME Prototype (LEAP) project,
funded by the U.S. Department of Energy, Office of Science, Office of Biological and Environmental Research,
under contract DE-AC52-06NA25396.
Los Alamos Report LA-UR-17-24044

\section*{References}


\begin{thebibliography}{10}
\expandafter\ifx\csname url\endcsname\relax
  \def\url#1{\texttt{#1}}\fi
\expandafter\ifx\csname urlprefix\endcsname\relax\def\urlprefix{URL }\fi
\expandafter\ifx\csname href\endcsname\relax
  \def\href#1#2{#2} \def\path#1{#1}\fi

\bibitem{Thuburn08}
J.~Thuburn, Some conservation issues for the dynamical cores of nwp and climate
  models, J. Comp. Phys. 227 (2008) 3715--3730.

\bibitem{PG17}
A.~Palha, M.~Gerritsma, A mass, energy, enstrophy and vorticity conserving
  ({MEEVC}) mimetic spectral element discretization for the 2{D} incompressible
  {N}avier-{S}tokes equations, J. Comp. Phys. 328 (2017) 200--220.

\bibitem{CS12}
C.~Cotter, J.~Shipton, Mixed finite elements for numerical weather prediction,
  J. Comp. Phys. 231 (2012) 7076--7091.

\bibitem{MC14}
A.~McRae, C.~Cotter, Energy- and enstrophy-conserving schemes for the
  shallow-water equations, based onmimetic finite elements, Q. J. R. Meteorol.
  Soc. 140 (2014) 2223--2234.

\bibitem{NSC16}
A.~Natale, J.~Shipton, C.~Cotter, Compatible finite element spaces for
  geophysical fluid dynamics, Dynamics and Statistics of the Climate System 1
  (2016) 1--31.

\bibitem{AL81}
A.~Arakawa, V.~Lamb, A potential enstrophy and energy conserving scheme for the
  shallow water equations, Mon. Wea. Rev. 109 (1981) 18--36.

\bibitem{Salmon04}
R.~Salmon, Poisson-bracket approach to the construction of energy- and
  potential-enstrophy-conserving algorithms for the shallow-water equations, J.
  Atmos. Sci. 61 (2004) 2016--2036.

\bibitem{Salmon07}
R.~Salmon, A general method for conserving energy and potential enstrophy in
  shallow-water models, J. Atmos. Sci. 64 (2007) 505--531.

\bibitem{CT14}
C.~Cotter, J.~Thuburn, A finite element exterior calculus framework for the
  rotating shallow-water equations, J. Comp. Phys. 257 (2014) 1506--1526.

\bibitem{TF10}
M.~Taylor, A.~Fournier, A compatible and conservative spectral element method
  on unstructured grids, J. Comp. Phys. 229 (2010) 5879--5895.

\bibitem{Gerritsma11}
M.~Gerritsma, {Edge Functions for Spectral Element Methods}, in: Spectral and
  High Order Methods for Partial Differential Equations, Vol.~76 of Lecture
  Notes in Computational Science and Engineering, Springer, 2011, pp. 199--207.

\bibitem{KG13}
J.~Kreeft, M.~Gerritsma, Mixed mimetic spectral element method for stokes flow:
  A pointwise divergence-free solution, J. Comp. Phys. 240 (2013) 284--309.

\bibitem{Vallis06}
G.~Vallis, Atmospheric and Oceanic Fluid Dynamics: Fundamentals and Large-Scale
  Circulation, Cambridge University Press, Cambridge, 2006.

\bibitem{TRSK09}
J.~Thuburn, T.~Ringler, W.~Skamarock, J.~Klemp, Numerical representation of
  geostrophic modes on arbitrarily structured {C}-grids, J. Comp. Phys. 228
  (2009) 8321--8335.

\bibitem{KPG11}
J.~Kreeft, A.~Palha, M.~Gerritsma, Mimetic framework on curvilinear
  quadrilaterals of arbitrary order, ArXiV 1111.4304, 2011.

\bibitem{Hiemstra14}
R.~Hiemstra, D.~Toshniwal, R.~Huijsmans, M.~Gerritsma, High order geometric
  methods with exact conservation properties, J. Comp. Phys. 257 (2014)
  1444--1471.

\bibitem{brezzi1991mixed}
F.~Brezzi, M.~Fortin, {Mixed and Hybrid Finite Element Methods}, Vol.~15 of
  Springer Series in Computational Mathematics, Springer, 1991.

\bibitem{BoffiMixedFiniteElements2013}
D.~Boffi, F.~Brezzi, M.~Fortin, {Mixed Finite Element Methods and
  Applications}, Vol.~1 of Springer Series in Computational Mathematics,
  Springer, 2013.

\bibitem{ArnoldBoffiBonizzoni}
D.~Arnold, D.~Boffi, F.~Bonizzoni, Finite element differential forms on
  curvilinear cubic meshes and their approximation properties, Numer. Math. 129
  (2015) 1--20.

\bibitem{ArnoldAwanou}
D.~Arnold, G.~Awanou, Finite element differential forms on cubical meshes,
  Mathematics of Computation, 83(288) (2013) 1551--1570.

\bibitem{ArnoldBoffiFalk}
D.~Arnold, D.~Boffi, R.~Falk, Quadrilateral h(div) finite elements, SIAM J.
  Numer. Anal. 43(6) (2005) 2429--2451.

\bibitem{bossavit_japan_computational_1}
A.~Bossavit, {Computational electromagnetism and geometry: (1) Network
  equations}, Journal of the Japan Society of Applied Electromagnetics 7~(2)
  (1999) 150--159.

\bibitem{bossavit_japan_computational_2}
A.~Bossavit, {Computational electromagnetism and geometry: (2) Network
  constitutive laws}, Journal of the Japan Society of Applied Electromagnetics
  7~(3) (1999) 294--301.

\bibitem{bossavit_japan_computational_3}
A.~Bossavit, {Computational electromagnetism and geometry: (3) Convergence},
  Journal of the Japan Society of Applied Electromagnetics 7~(4) (1999)
  401--408.

\bibitem{bossavit_japan_computational_4}
A.~Bossavit, {Computational electromagnetism and geometry: (4) From degrees of
  freedom to fields}, Journal of the Japan Society of Applied Electromagnetics
  8~(1) (2000) 102--109.

\bibitem{bossavit_japan_computational_5}
A.~Bossavit, {Computational electromagnetism and geometry: (5) The ``Galerkin
  Hodge''}, Journal of the Japan Society of Applied Electromagnetics 8~(2)
  (2000) 203--209.

\bibitem{arnold2006finite}
D.~N. Arnold, R.~S. Falk, R.~Winther, {Finite element exterior calculus,
  homological techniques, and applications}, Acta Numerica 15 (2006) 1--155.

\bibitem{arnold2010finite}
D.~N. Arnold, R.~S. Falk, R.~Winther, {Finite element exterior calculus: from
  Hodge theory to numerical stability}, Bulletin of the American Mathematical
  Society 47~(2) (2010) 281--354.

\bibitem{Palha2014}
A.~Palha, P.~Rebelo, R.~Hiemstra, J.~Kreeft, M.~Gerritsma, {Physics-compatible
  discretization techniques on single and dual grids, with application to the
  Poisson equation of volume forms}, Journal of Computational Physics 257
  (2014) 1394--1422.

\bibitem{robidoux-polynomial}
N.~Robidoux, {Polynomial histopolation, superconvergent degrees of freedom, and
  pseudospectral discrete Hodge operators}, Unpublished:
  http://people.math.sfu.ca/$\sim$nrobidou/public\_html/prints/histogram/histogram.pdf.

\bibitem{abraham_diff_geom}
R.~Abraham, J.~E. Marsden, T.~Ratiu, {Manifolds, Tensor Analysis, and
  Applications}, Vol.~75 of Applied Mathematical Sciences, Springer, 2001.

\bibitem{frankel}
T.~Frankel, {The Geometry of Physics}, 2nd Edition, Cambridge University Press,
  2004.

\bibitem{SB85}
R.~Sadourny, C.~Basdevant, Parameterization of subgrid scale barotropic and
  baroclinic eddies in quasi-geostrophic models: Anticipated potential
  vorticity method, J. Atmos. Sci. 42 (1985) 1353--1363.

\bibitem{ER17}
C.~Eldred, D.~Randall, Total energy and potential enstrophy conserving schemes
  for the shallow water equations using hamiltonian methods – part 1:
  Derivation and properties, Geosci. Model Dev. 10 (2017) 791--810.

\bibitem{KS05}
G.~Karniadakis, S.~Sherwin, {Spectral/hp Element Methods for Computational
  Fluid Dynamics, Second Edition}, Oxford, 2005.

\bibitem{MST12}
T.~Melvin, A.~Staniforth, J.~Thuburn, Dispersion analysis of the spectral
  element method, Q. J. R. Meteorol. Soc. 138 (2012) 1934--1947.

\end{thebibliography}

\end{document}